\documentstyle[
10pt, epsf,stmaryrd,amscd,amstex,amssymb]{amsart}
\newtheorem{Thm}{Theorem}
\newtheorem{Cor}{Corollary}
\newtheorem{Lem}{Lemma}
\newtheorem{Prop}{Proposition}

\theoremstyle{remark}
\newtheorem{Rem}{Remark}
\newtheorem{Def}{Definition}
\newtheorem{Ex}{Example}

\newcommand{\Ql}{{\overline {{\Bbb Q}_l} }}

\newcommand{\Fl}{{{\cal F}\ell}}

\newcommand{\bu}{\bullet}

\newcommand{\To}{\longrightarrow}
\newcommand{\iso}{{\widetilde \longrightarrow}}
\newcommand{\isol}{{\widetilde \longleftarrow}}
\newcommand{\imbed}{\hookrightarrow}

\newcommand{\<}{\langle}
\renewcommand{\>}{\rangle}

\newcommand{\bl}{\ll} 
\newcommand{\br}{\gg} 

\def\square{\hbox{\vrule\vbox{\hrule\phantom{o}\hrule}\vrule}}
\newcommand{\epf}{\square}

\renewcommand{\P}{{\cal P}}

\newcommand{\Sch}{{\mathcal S}}

\newcommand{\LL}{{\mathcal L}}

\newcommand{\PCoh}{{\cal P}Coh}

\newcommand{\PIW}{\cal P_{{\cal {IW}}}}

\newcommand{\Phitil}{\tilde{\Phi}}
\newcommand{\Phit}{\tilde{\Phi}}

\newcommand{\N}{{\cal N}}
\newcommand{\Ntil}{{\tilde{\cal N}}}

\newcommand{\Nt}{{\tilde{\cal N}}}
\newcommand{\ICtil}{\widetilde{IC}}

\newcommand{\Ft}{{\tilde F}}

\newcommand{\Ht}{{\widetilde{H}}}

\newcommand{\pit}{{\widetilde{\pi}}}

\renewcommand{\j}{{\frak j}}

\newcommand{\oplusl}{\bigoplus\limits}
\newcommand{\cupl}{\bigcup\limits}

\renewcommand{\b}{{\frak b}}

\newcommand{\bq}{{\bf q}}

\newcommand{\g}{{\frak g}}
\newcommand{\n}{{\frak n}}

\newcommand{\LG}{{ G\check{\ }}}
\newcommand{\LB}{{B\check{\ }}}

\newcommand{\LN}{{N\check{\ }}}

\renewcommand{\O}{{\cal O}}

\newcommand{\F}{{\cal F}}
\newcommand{\G}{{\cal G}}
\newcommand{\A}{{\cal A}}
\newcommand{\cB}{{\cal B}}

\newcommand{\C}{{\cal C}}
\newcommand{\E}{{\cal E}}
\newcommand{\V}{{\cal V}}
\newcommand{\cW}{{\cal W}}

\newcommand{\cons}{{\overline{\underline{{\Bbb Q}_l}}}}

\newcommand{\Zet}{{\Bbb Z}}
\newcommand{\Ce}{{\Bbb C}}

\renewcommand{\c}{{\underline {c}}}

\newcommand{\bI}{{\bf I}}

\newcommand{\LGO}{{\bf \LG _O}}

\newcommand{\LGK}{{\bf \LG_{\Fqbar((t))}}}

\newcommand{\Db}{D^b}

\newcommand{\Pone}{{\mathbb P}^1}
\newcommand{\Aone}{{\mathbb A}^1}

\newcommand{\Gm}{{\mathbb G}_m}

\newcommand{\Fq}{{\mathbb F}_q}
\newcommand{\Fqbar}{{\bar{\mathbb F}_q}}

\newcommand{\U}{{\mathsf U}}
\newcommand{\Uq}{{\U_q}}
\newcommand{\Bq}{{{\mathsf B}_q}}
\newcommand{\Uqmod}{\U_q{\textrm -mod^0}}
\newcommand{\Weyl}{{\rm Weyl}}
\newcommand{\coWeyl}{{\rm coWeyl}}

\newcommand{\lTo}{\longleftarrow}

\newcommand{\D}{{\mathfrak D}}

\renewcommand{\i}{{\mathfrak i}}

\renewcommand{\k}{{\bf k}}

\newcommand{\Hom}{{\rm Hom}}
\newcommand{\End}{{\rm End}}
\newcommand{\Ext}{{\rm Ext}}
\newcommand{\Ob}{{\rm Ob}}

\newcommand{\nab}{\nabla}

\newcommand{\del}{\Delta}

\newcommand{\codim}{{\rm co}\dim}

\newcommand{\la}{\lambda}

\newcommand{\La}{\Lambda}

\newcommand{\al}{\alpha}

\newcommand{\pr}{\pi}

\newcommand{\proofpt}{{\it Proof. }}

\author{Roman Bezrukavnikov}
\title
[Cohomology of tilting modules and $t$-structures]
{Cohomology of tilting modules over quantum groups and
$t$-structures on derived categories of coherent sheaves}
\begin{document}
\maketitle
\centerline{\em To the memory of my father.}
\medskip



\begin{abstract}{\footnotesize}
\noindent The paper is  concerned with   cohomology of the small
quantum group at a root of unity, and of its upper triangular
subalgebra, with coefficients in a tilting module. It turns out to
be related to irreducible objects in the heart of a certain
$t$-structure on the derived category of equivariant coherent
sheaves on the Springer resolution, and to equivariant coherent IC
sheaves on the nil-cone. The support of the cohomology is described
in terms of cells in affine Weyl groups. The basis in the
Grothendieck group provided by the cohomology modules is shown  to
coincide with  the Kazhdan-Lusztig basis,
 as predicted by J.~Humphreys and V.~Ostrik.

The proof is based on the results of \cite{ABG}, \cite{AB} and \cite{B},
which allow us to reduce the question to purity of IC sheaves
on affine flag varieties.
\end{abstract}

\tableofcontents

\section{Introduction.}
\subsection{Quantum groups, tilting modules and cohomology.}
Let $G$ be a semi-simple complex algebraic group of adjoint type,
$N\subset B\subset G$ be  a maximal unipotent and a Borel subgroups,
and $\n\subset \b\subset \g$ be the corresponding  Lie algebras. Let
$\N\subset \g$ be the nilpotent cone, let $\Nt=T^*(G/B)=G\times _B\n
=\{(gB, x)\ |\ x\in Ad(g)(\n)\}$ and let $\pi:\Nt\to \N$
 be the Springer-Grothendieck map.
For an algebraic group $H$ acting on an algebraic  variety $X$ we
let $Coh^H(X)$ denote the category of $H$-equivariant coherent
sheaves on $X$, and $D^H(X)=D^b(Coh^H(X))$ be the bounded derived
category.


Let $q\in \Ce$ be a primitive root of unity
of order $l$. We assume that $l$ is odd, and prime to 3 if $\g$ has a factor
of type $G_2$; and also that $l$ is greater than the Coxeter number.
 Let $\Uq$ be the Lusztig quantum enveloping algebra
\cite{Lbook}. We have finite dimensional subalgebras $b_q\subset u_q
\subset \Uq$, where $u_q$ is the so-called ``small quantum group'',
and $b_q\subset u_q$ is the upper triangular subalgebra.

 For an augmented
$\Ce$-algebra $A$ we write $H^\bu(A)=\Ext_A(\Ce,\Ce)$, and
$H^\bu(A,M)=\Ext_A(\Ce, M)$ for an $A$ module $M$; thus $H^\bu(A)$
 is a graded
 associative algebra, and $H^\bu(A,M)$ is a graded
module over this algebra.

We consider the category $\Uq{\textrm -mod}$ of finite dimensional
graded $\Uq$ modules; we let $\Uqmod\subset \Uq{\textrm -mod}$ be
the block of the trivial module, see, e.g., \cite{ABG}, \S 3.4 (cf.
also \cite{AG}, \S 1.2).

 According to \cite{GK}
we have $H^\bu(u_q)\cong \O(\N)$.
Thus for a $u_q$ module $M$ we get an $H^\bu(u_q)=\O(\N)$ module
 $H(M)=H^\bu(u_q,M)$.
Moreover, it is explained in \cite{GK} that if $M$ is an $\Uq$
module then the graded $\O(\N)$ module $H(M)$ is $G$-equivariant,
thus we get a functor $H:\Uq{\textrm -mod}\to Coh^{G\times
\Gm}(\N)$, $H(M)=H^\bu(M)$, where the multiplicative group $\Gm$
acts on $\N$ by
$t:x\mapsto t^2x$. 

In this paper we provide a reasonably explicit description of $H(M)$
when $M$ is a {\it tilting} object of $\Uqmod$. To state this
description we have to recall the category $\PCoh$ of
$G$-equivariant coherent perverse sheaves on $\N$ with respect to
the middle perversity, see \cite{izvrat}. It is the heart of a
certain $t$-structure on $D^G(\N)$. For each pair $(O,\LL)$ where
$O\subset \N$ is a $G$ orbit, and $\LL$ is an irreducible
$G$-equivariant vector bundle on $O$ there is a unique irreducible
object $IC_{O,\LL}$ supported on the closure of $O$, and satisfying
$IC_{O,\LL}|_{O}=\LL[-\frac{\codim O}{2}]$; these are all of the
irreducible objects of $\PCoh$. We will prove the following

\begin{Thm}\label{cohoti}
 Let $T$ be an indecomposable tilting object of $\Uqmod$.
Then either $H(T)=0$, or there exists a(n obviously unique) pair
$(O_T,\LL_T)$
 as above, such that we have an isomorphism of $G$-equivariant $\O(\N)$
modules $H(T)\cong \Gamma (IC_{O_T,\LL_T})$; here $\Gamma=\oplusl_i
R^i\Gamma $ denotes the total cohomology (derived global sections)
functor.
\end{Thm}

See  Corollary \ref{HTlambda} in section \ref{sect32}
for a more precise statement. The
latter  implies, via a result of \cite{B}, a conjecture of
J.~Humphreys (proved for type $A$ in Ostrik's thesis, see \cite{Ohu}), which
describes the
 support of $H(T)$ in terms of 2-sided cells
in the affine Weyl group (Corollary \ref{Humphreys} in section
\ref{sect32}). Our results also yield some conjectures by V.~Ostrik,
see Remark \ref{proO}. These applications have partly motivated the
present work.

\bigskip

In fact, the above statements will be derived from
 the following stronger results.

First of all, besides of the objects $H(T)\in Coh^G(\N)$ we will
also describe (though somewhat less explicitly) the objects $\Ht
(T)\in Coh^G(\Nt)$ carrying more information. These are defined as
follows. By \cite{GK} we have $H^\bu(b_q)=\O(\n)=Sym(\n^*)$; thus
for a $b_q$ module $M$ we get an $\O(\n)$ module $H^\bu(b_q,M)$. If
$M$ is equipped with an action of the Borel subalgebra $\Bq\subset
\Uq$ then the $\O(\n)$ module  $H^\bu(b_q,M)$ is equivariant with
respect to the adjoint action of the Lie algebra $\b$. If $M$ is one
dimensional, an explicit calculation (see \cite{GK}) shows that this
module is locally finite; it follows that the same is true for any
complex $M$ of $\Bq$ modules whose cohomology is finite dimensional.
Thus for such $M$ the $\n$ action on $H^\bu(b_q,M)$  integrates to
an action of the algebraic group $N$; furthermore, if $M$ admits a
grading by weights compatible with the $\Bq$ action, then the action
of $\b$ integrates to an action of $B$,
 so we get a functor from the bounded derived category of finite
dimensional graded $\Bq$ modules to $Coh^B(\n)$. Composing it with
the restriction functor from $\Uq$ to $\Bq$ we get a functor
$\Ht:\Uq{\textrm -mod} \to Coh^B(\n)\cong Coh^G(\Nt)$, where the
last equivalence is the induction functor, inverse to $\F\mapsto
\F|_\n$ (recall that $\Nt=G\times_B \n$).

The description of $\Ht(T)$ for a tilting object $T\in \Uqmod$ is as
follows.  We define a $t$-structure on the bounded derived
categories $D^{G\times \Gm}(\Nt)$,
 $D^{G}(\Nt)$,
which we call the {\it exotic $t$-structure}; we call objects of its
heart {\it exotic sheaves}. We
  prove that the sheaves
$\tilde H(T)$, where $T$ is an indecomposable tilting object of $\Uqmod$,
are precisely the cohomology sheaves of irreducible exotic sheaves.
 It turns out that for an irreducible exotic sheaf
$E$ the object $R\pi_*(E)$ is either zero, or isomorphic to
$IC_{O,\LL}$ for some $(O,\LL)$.
 This allows one to deduce Theorem \ref{cohoti}.

Actually, all the above results about cohomology of $\Uq$ modules
will be deduced from  stronger statements involving isomorphisms of objects
in the derived category of coherent sheaves. To state the latter we recall
that (the first part of) \cite{ABG} provides a
triangulated functor $\Psi:D^{G\times \Gm}(\Nt)\to D^b(\Uqmod)$
which is ``almost an equivalence'' (see \eqref{equiv}
 below for a precise statement), here $\Gm$ acts on $\Nt$ by
$t:(gB,x)\mapsto (gB, t^2x)$.

The following statement is the  central point of the present paper, see
Theorem \ref{ET} below for a more precise statement.

\begin{Thm}\label{gogo} For any irreducible exotic
sheaf $E\in D^{G\times \Gm}(\Nt)$ we have $\Psi(E)\cong T[n]$ for
some indecomposable tilting object $T\in \Uqmod$ and $n\in \Zet$.
Every indecomposable tilting $T\in \Uqmod$ has the form $\Psi(E)$
for a unique irreducible exotic sheaf
 $E\in  D^{G\times \Gm}(\Nt)$.
\end{Thm}

This will be deduced from purity of irreducible perverse sheaves
on the affine flag variety for the Langlands dual group via results
of \cite{AB}.

\subsection{Koszul duality conjecture.}\label{koszul_sect}
 In this section we attempt
to explain the origin of our method.

To clarify the idea behind the argument we gather together
 various functors, some
of which are
used in this paper:

$$D^b(\P_{sph})\underset{(I)}{\cong}
D^b(\Uqmod)\underset{(II)}{\cong} DGCoh^G(\Nt)
\overset{(III)}{\lTo}
 D^{G\times \Gm}(\Nt)\overset{(IV)}{\To} D^G(\Nt)
\underset{(V)}{\cong}   D^b(\P_{asph}) $$

Here $D^b(\P_{sph})$ (where ``sph'' stands for ``spherical'') is the
category of perverse sheaves on the affine Grassmannian for the
Langlands dual group $\LG$ smooth along the Schubert stratification;
the  equivalence (I) is defined in \cite{ABG}. Unlike the other
functors appearing above, this equivalence comes from an {\it
equivalence of abelian categories.} Neither (I), nor the category
$\P_{sph}$ are used in the present paper, they are
 mentioned here to clarify the picture; thus we will use only the first,
``algebraic'' part of \cite{ABG}.

The  equivalence (II) is also introduced in \cite{ABG}. The category
$DGCoh^G(\Nt)$ here is the (derived) category of DG modules over a
certain DG  algebra with zero differential, such that the category
of modules over its cohomology algebra is naturally identified with
$Coh^G(\Nt)$. We do not go into details here, referring the
interested reader to \cite{ABG}. Instead we use the composition of
functors (III) and (II). We remark that providing the composite
functor
 $D^{G\times \Gm}(\Nt)\to
 D^b(\Uqmod)$ satisfying properties \eqref{twist}, \eqref{equiv}
below is essentially equivalent to providing the equivalence (II).

(IV) is the forgetful functor, restricting the equivariance from $G\times
\Gm$ to $G$.

Finally, $D^b(\P_{asph})$ (where ``asph'' stands for
``anti-spherical'') is a certain category of perverse sheaves on
the affine flag variety of the Langlands dual group $\LG$, and (V)
is proved in \cite{AB} (see section \ref{onAB} below for more
details). This equivalence will be used in order to deduce the key
Positivity Lemma (Lemma \ref{posi_lem}
in section \ref{posisect}) below from purity for perverse sheaves on the affine
flag variety.

\medskip

The shortest (though, perhaps, the least elementary) description of
the  exotic $t$-structure on $D^G(\Nt)$ is as follows: it is the
transport of the tautological $t$-structure on $D^b(\P_{asph})$ by
means of the equivalence (V); the exotic $t$-structure on
$D^{G\times \Gm}(\Nt)$ is then the unique $t$-structure compatible
with the exotic $t$-structure on $D^G(\Nt)$. Thus an irreducible
exotic sheaf
 goes to an irreducible object of
$D^b(\P_{asph})$ under the composition (V)$\circ$(IV). Comparing
this with Theorem \ref{gogo} we discover a relation between tilting
objects in $\P_{asph}$ and irreducible objects in $\P_{sph}$.
 This motivates the following conjecture, which inspired
the methods of the present work.

 To state  the conjecture we need
the following concept. We will say that two abelian categories $\A$
and $\cB$ are {\em derived Koszul equivalent} if there exist graded
versions\footnote{See, e.g., \cite{BGS}, \S 4.3 for the definition.}
$\A^{gr}$, $\cB^{gr}$ and an equivalence of derived categories
$\kappa:\Db(\A^{gr})\cong \Db(\cB^{gr})$, which satisfies
$\kappa(M(1))\cong \kappa(M)(1)[1]$, where $M\mapsto M(1)$ is the
shift of grading functor. The interested reader is referred to
\cite{BGS} for examples and comments.

\medskip

{\bf Conjecture.} \emph{The category of perverse sheaves on the
affine flag variety of a semi-simple group, which are smooth along
the Schubert stratification is derived Koszul equivalent to itself;
the corresponding equivalence interchanges irreducible and tilting
objects. There exists also a parabolic -- singular version (in the
sense of \cite{BGS}) of this duality, whose particular case is a
derived Koszul equivalence between $\P_{sph}$ and $\P_{asph}$; the
corresponding equivalence interchanges tilting and irreducible
objects.}

\medskip

This Conjecture is an affine analogue of the main results of
\cite{S}, \cite{BGS}, or rather of a variant of the latter provided
by \cite{BG}.

In fact,  existence of the derived Koszul equivalence between
$\P_{sph}$ and $\P_{asph}$ sending irreducibles to tiltings  follows
from the above mentioned results of \cite{AB}, \cite{ABG}, and the
present paper. It is also plausible that the first part of the
conjecture can be proven by the methods of \cite{S} complemented by
the geometric counterpart of translation functors. The latter is a
certain collection of exact endofunctors of the category of
unipotently monodromic perverse sheaves on the basic affine space of
a Kac-Moody group, which in the case of a finite dimensional group
reduces to the usual reflection functors. They can be defined as
convolution with tilting sheaves (cf. \cite{BG}, Theorem 6.10); some
other constructions will appear in a forthcoming paper by
I.~Mirkovi\' c.

\medskip

Finally, let us mention that the exotic $t$-structure on the derived
category of coherent sheaves on $\Nt$ appears in some other
contexts, e.g., it is related to the $t$-structures appearing in
\cite{BMR}, and, conjecturally,\footnote{The conjecture stems from
discussions with D.~Gaitsgory.} to modules over the affine Lie
algebra on the critical level studied by Frenkel and Gaitsgory, cf.
\cite{FG}.

\subsection{Plan of the paper.} In section \ref{coh}, after some
 preliminaries about mutations of exceptional sets and $t$-structures
defined by  generating exceptional sets, we define the exotic $t$-structure
on $D^G(\Nt)$, $D^{G\times \Gm}(\Nt)$; we also relate it to the
middle perversity $t$-structure on $D^G(\N)$. We finish the section
by stating the Positivity Lemma.
 In section \ref{qua}, modules
over the quantum group at a root of unity appear: we recall the results
of \cite{ABG} which relate their derived category to $D^{G\times \Gm}(\Nt)$.
We then use the Positivity Lemma to show that a functor constructed in
\cite{ABG} sends irreducible exotic sheaves to tilting modules over a quantum
group (up to a homological shift). In section \ref{onAB} we recall from
\cite{AB} the
equivalence between $D^G(\Nt)$ and a certain derived category of perverse
sheaves on the affine flag variety for the Langlands dual group.
It sends irreducible exotic sheaves to irreducible perverse sheaves,
and the grading on $\Hom$ spaces
 coming from the $\Gm$-equivariant structure is related
to Frobenius weights. This allows to deduce Positivity Lemma from
purity of irreducible perverse sheaves (see also Remark \ref{purpur}
in section \ref{posisect}).

\subsection{Notations.}\label{nota}
 Fix an algebraically closed ground
field $\k$ of characteristic zero.

\subsubsection{Notations related to $G$}\label{notaG}
$G$ is a semi-simple algebraic group; in section \ref{qua} it will
be assumed to be adjoint. The Springer map $\pi:\Nt\to \N$ was
recalled above.

 Let $\Lambda$ be the weight lattice
of $G$,  $\Lambda^+\subset \Lambda$ be the set of dominant weights
and $R^+\subset \La$ be the subsemigroup (with zero) generated by
positive roots. We have the standard partial order on $\La$,
$\la\preceq \mu$ if $\mu-\la\in R^+$.

 For $\lambda\in \Lambda^+$ we let $V_\lambda$ be the
corresponding irreducible representation of $G$ or of the Lie
algebra $\g$. We let $W_f$ be the Weyl group of $G$, and
$W=W_f\ltimes \Lambda$ be the extended affine Weyl group. Let
$\ell:W\to \Zet_{\geq 0}$ be the length function. For $\lambda\in
\Lambda$ let $\delta_\lambda$ denote the minimal length of an
element $w\in W_f$ such that $w(\lambda)\in \Lambda^+$.


\subsubsection{Notations for triangulated categories}\label{notatr}
For a triangulated category $\D$ and $X,Y\in \D$ we will use the
notation $\Hom^n(X,Y)=\Hom(X,Y[n])$, $\Hom^\bu(X,Y)=\oplusl_{n\in \Zet}
\Hom(X,Y[n])$.

For a set $S$ of objects in $\D$
we will denote by $\bl S \br$ the full triangulated subcategory
generated by $S$; and by $\< S\>$ the full subcategory generated by $S$ under
 extensions. Thus $\bl S\br$ is the smallest strictly full triangulated
subcategory of $\D$ containing $S$; and $\< S\>$ is the smallest
 strictly full subcategory containing $S$ and closed under extensions
(we say that $\C\subset \D$ is closed under extensions if for any $X,Y\in \C$
and any exact triangle $X\to Z\to Y\to X[1]$ we have $Z\in \D$).

If $\C\subset \D$ is a full triangulated subcategory we can form
the quotient triangulated category $\D/\C$; we will write $X\cong Y\mod \C$
for $X,Y\in \D$ meaning that the images of $X$ and $Y$ in $\D/\C$
are isomorphic.

 For a category $\C$ let $[\C]$ be the set
of isomorphism classes of objects in $\C$.
 If $S_1,S_2$ are
subsets of  $[\D]$, then $S_1*S_2$ denotes the subset of
$[\D]$ consisting of
 classes of all objects $Z$, for which there exists an
exact triangle $X\to Z\to Y\to X[1]$ with $[X]\in S_1$, $[Y]\in
S_2$.

\medskip

Direct and inverse image functors on the derived categories of (quasi)coherent
sheaves etc. are understood to be the corresponding derived functors,
unless stated otherwise.

\subsection{Acknowledgements.}
This is an outgrowth of the project initiated during the IAS
special year in Representation Theory (1998/99), as a result of
various conversations with several mathematicians, especially with
M.~Finkelberg and I.~Mirkovi\' c; I am much indebted to them. I also
thank D.~Gaitsgory, V.~Ginzburg and V.~Ostrik for inspiring communications.
I thank J.~Humphreys, J.~Kamnitzer, I.~Mirkovi\' c and
 the referees for pointing out many typos and exposition lapses.

The author was supported by the NSF grant DMS-0071967, and was
employed by the Clay Mathematical Institute at the time when ideas
of this paper were being worked out; he was partially
 supported by DARPA grant
HR0011-04-1-0031 during the final stage of his work on the paper.

\section{Coherent sheaves: exotic $t$-structure, and Positivity Lemma.}
\label{coh} In this section we define the exotic $t$-structures on
$D^G(\Nt)$, $D^{G\times \Gm}(\Nt)$.  The construction will use the
notion of an exceptional set in a triangulated category; we start by
recalling this notion.

\subsection{Exceptional sets and
quasi-hereditary hearts.}
Most of the material in this section is borrowed from \cite{BK} or
\cite{BGS}.

\subsubsection{Admissible subcategories}
For a subcategory
$\C$ in an additive category $\D$ let us (following \cite{BK}) write
 $\C^\perp=(\C^\perp)_\D$ (respectively $^\perp \C=(^\perp \C)_\D$)
for the strictly full subcategory in $\D$ consisting of objects $X$ for which
$\Hom(A,X)=0$ (respectively $\Hom(X,A)=0$) for all $A\in \C$. The subcategories
$\C^\perp$, $^\perp \C$ are called respectively the right and left orthogonal
of $\C$.

A full triangulated subcategory $\C\subset \D$ is called \emph{right}
(respectively, \emph{left}) \emph{admissible} if the following equivalent
conditions hold (see, e.g., \cite{BK}, Proposition 1.5).

a) The inclusion functor $\C\imbed \D$ has a right (respectively, left)
adjoint.

b) $\D=\<\C, \C^\perp\>$ (respectively, $\D=\< \C, ^\perp \C\>$).

c)  $[\D]=[\C]* [\C^\perp]$ (respectively, $[\D]=[^\perp \C]* [\C] $).

(See \ref{notatr} for notations.)
\medskip

If $\C\subset \D$ is a right (respectively, left) admissible subcategory,
then $\C^\perp$ (respectively, $^\perp\C$) is left (right) admissible.

\subsubsection{Exceptional sets}
Let $\D$ be a $\k$-linear triangulated category. We assume it is {\it of finite
type}, i.e. the $\k$ vector space $\Hom^\bu(X,Y)$ is finite dimensional for any
objects $X,Y$ of $\D$.

An ordered subset  $\nab=\{ \nab^i,\;  i\in I\}$ of $\Ob(\D)$ is
called {\it exceptional} if we have $\Hom^\bu (\nab^i,
\nab^j)=0$ for $i<j$;
 $\Hom ^{n}(\nab^i,\nab^i)=0$ for $n\ne 0$,
 and $\End(\nab^i)=\k$.

Let $\nab= \{ \nab^i,\; i\in I\}$ be an exceptional set.
For $i\in I$ let
 $\D_{<i}$ 
 be the full triangulated
subcategory generated by  $\nab_j$, $j<i$.

Let  $\del=\{ \del_i,\;  i\in I\}$ be another subset of
$\Ob(\D)$ (in bijection with $\nab$).
We say that $\del$ is  {\it dual} to $\nab$ if
 \begin{equation}\label{YX}
\Hom^\bu (\del_n, \nab^i)=0 \ \ {\rm for}\ \ n>i;
\end{equation}
and there exists an isomorphism
 \begin{equation}\label{XY}
\del_n \cong \nab^n \mod \D_{<n}.
\end{equation}

It is easy to see  then that  $\Hom^\bu(\del_i,\nab^i)=\k$, and
$\Hom^\bu(\del_i,\nab^j)=0$ for $i\ne j$; that the dual set
equipped with the opposite ordering is exceptional; and that it is
unique if it exists (up to an isomorphism, which is fixed uniquely by
fixing  \eqref{XY});
 see, e.g., Lemma 2 in \cite{nilpokon}. Moreover, we have

\begin{Prop}\label{exist_dual}
Let $\nab\subset Ob(\D)$ be a finite exceptional set.

a)
The triangulated
 subcategory $\C=\bl \nab\br$ generated by $\nab$ is both left and right
admissible.

b) The dual exceptional set exists.
\end{Prop}

\proofpt (a) is  Corollary 2.10 in \cite{BK}. (b) follows from (a):
set $\C_n =^\perp \< \nab^1, \dots, \nab^{n-1}\>$; then $\C_n$
is a right admissible subcategory, thus there exists a functor
$\Pi_n:\D\to \C_n$
right adjoint to the inclusion. 
 We set $\del_n=\Pi_n(\nab^n)$;
the desired properties of $\del_n$ follow by a straightforward verification.
\epf

\subsubsection{$t$-structure of an exceptional set and quasi-hereditary
categories}\label{qhc}  Let $\nab=(\nab^i)$ be an exceptional set in
a finite type triangulated $\k$-linear category $\D$; here $i$ runs
over $\Zet_{>0}$, or $[1,..,n]$. Let $\del_i$ be the dual
exceptional set. Assume that $\D=\bl \nab \br$.

We refer, e.g., to \cite{nilpokon} for a proof of the next
Proposition; here we remark only that the $t$-structure is
obtained  by applying gluing of $t$-structures construction from
\cite{BBD}, \S 1.4.

\begin{Prop} \label{propexcep}
a) There exists a unique $t$-structure $(\D^{\geq 0}, \D^{<0})$ on
$\D$, such that $\nab^i\in \D^{\geq 0}$; $\del^i\in \D^{\leq 0}$.
Moreover, $\D^{\geq 0}$, $\D^{<0}$ are given by
\begin{equation}\label{defDbo}
  \D^{\geq 0}=\< \{\nab^i[d]\ ,\ i\in I, \; d\leq 0\} \>;
\end{equation}
\begin{equation}\label{defDme}
\D^{< 0}=\< \{\del_i[d]\ ,\ i\in I, \; d> 0\} \>.
\end{equation}

b) The $t$-structure is bounded.

c) For $X\in \D$ we have $X\in  \D^{\geq 0}\iff \Hom^{<0}(\del_i,
X)=0\ \forall i$; 

$X\in D^{<0}\iff \Hom^{\leq 0 }(X, \nab^i)=0\ \forall i$.

d) Let $\A$ denote the heart of the above $t$-structure. Then every
object of $\A$ has finite length. For every $i$ the image $L_i$ of
the canonical arrow $\tau_{\geq 0}(\del_i)\to \tau_{\leq 0}(\nab^i)$
is irreducible. The objects $L_i$ are pairwise non-isomorphic, and
every irreducible object in $\A$ is isomorphic to $L_i$ for some
$i$. The $t$-structure induces a $t$-structure on $\D_i=\bl \{
\nab_1,\dots, \nab_i\} \br$; the heart of the latter, $\A_i$, is the
Serre subcategory in $\A$ generated by $L_1,\dots L_i$. The map
$\tau_{\geq 0}(\del_i)\to L_i$ is a projective cover of $L_i$ in
$\A_i$, and $L_i\to \tau_{\leq 0}(\nab^i)$ is an injective hull of
$L_i$ in $\A_i$. \epf

\end{Prop}

\begin{Rem}
An abelian category satisfying the properties summarized in Proposition
\ref{propexcep}(d) is called
a quasi-hereditary category (or Kazhdan-Lusztig type category, or highest
 weight category) see, e.g., \cite{PS} or \cite{BGS}.
\end{Rem}

\begin{Rem}\label{izvrat_puch}
The reader can keep in mind the following example.
Let $\D$ be a full subcategory in the bounded  derived category of
 sheaves of vector spaces
 on a reasonable topological space (or of the  ``derived category''
of $l$-adic
etale sheaves on a reasonable
scheme),
consisting of complexes whose cohomology is
smooth along a fixed (reasonable) stratification. Assume for simplicity that
the strata $\Sigma_i$ are connected and simply-connected, and
satisfy $H^{>0}(\Sigma_j)=0$; we
write
$j<i$ if $\Sigma_j$ lies in the closure of $\Sigma_i$.
 Let
$\j_i$  denote the embedding of  $\Sigma_i$ in the space.
Let $p_i$ be
arbitrary integers.
Then  objects $\nab^i=\j_*(\cons [p_i])$ form an exceptional set
generating
$\D$, and $\del_i=\j_!(\cons [p_i])$ is the dual set.
The $t$-structure of this exceptional set
is the perverse $t$-structure \cite{BBD} corresponding to
perversity $p=(p_i)$.
\end{Rem}

\subsubsection{Mutation of an exceptional set} Let $(I,\preceq)$
be an ordered set, and $\nab^i\in \D$, $i\in I$ be an exceptional
set. Let $\leq$ be another order on $I$; we assume that either $I$
is finite, or the ordered set $(I,\leq)$ is isomorphic to
$\Zet_{>0}$. Recall that $\D_{\leq i} =\bl\{ \nab^j\ |\ j\leq i\}
\br$; $\D_{< i} =\bl \{ \nab^j\ |\ j< i\} \br$.

Let $(I,\preceq)$ be an ordered set, and $\nab^i\in \D$, $i\in I$ be
an exceptional set. Let $\leq$ be another order on $I$; we assume
that either $I$ is finite, or the ordered set $(I,\leq)$ is
isomorphic to $\Zet_{>0}$. Recall that $\D_{\leq i} =\bl\{ \nab^j\
|\ j\leq i\} \br$; $\D_{< i} =\bl \{ \nab^j\ |\ j< i\} \br$.

\begin{Lem}\label{mutalem}
a) For $i\in I$ there exists a unique (up to a unique isomorphism)
object $\nab^i_{mut}$ such that $\nab^i_{mut}\in \D_{\leq i}\cap
\D_{<i}^\perp$, and $\nab^i_{mut}\cong \nab^i\mod \D_{<i}$.

b) The objects $\nab^i_{mut}$ form an exceptional set indexed by
$(I,\leq)$.

c) We have $\D_{\leq i} =\bl \{ \nab^j_{mut}\ |\ j\leq i\} \br$,
$\bl \nab\br=\bl \nab_{mut} \br$.
\end{Lem}

\proofpt The proof of uniqueness in (a) is standard. Let $\Pi_i$
denote the projection functor $\D_{\leq i}\to \D_{\leq i}/\D_{<i}$.
Let $\Pi^r_i$ denote the right adjoint functor, it exists by
Proposition \ref{exist_dual}(a).
 We set
$\nab^i_{mut}=\Pi_i^r \circ \Pi_i(\nab^i)$; it is immediate to check
that $\nab^i_{mut}$ satisfies the requirements of part (a) and forms
an exceptional set. Part (c) is then easily proved by induction.
\epf

\medskip

 We will say that the exceptional set $(\nab^i_{mut})$ is
 the $ \leq$ mutation of $(\nab^i)$. It is clear from the definition
 that it depends only on the second order $\leq$, but not on the
 original order $\preceq$.

\begin{Rem} The above notion of mutation of an exceptional set is
related to the action of the braid group on the set of exceptional
sets in a given triangulated category constructed in \cite{BK}
(this action is also called the action  by mutations).
Assume for simplicity that an exceptional set $\nab$ is finite. Then
a pair of orders $(\preceq, \leq)$ on $\nab$ defines a permutation
$\sigma$ of $\nab$. The $(\leq)$ mutation of $\nab$ can also be
obtained by the action of the element $\tilde \sigma$ of the braid
group on $\nab$; here $\tilde \sigma$ is the minimal length lifting
of $\sigma$ to an element of the braid group.
\end{Rem}

\subsubsection{Graded version}\label{grv}
Here we extend (in a rather straightforward manner) the above
results about exceptional sets to some infinite exceptional sets.
An example of the situation discussed in this subsection arises
from a graded quasi-hereditary algebra: the bounded derived
category of finitely generated modules is generated by a finite
exceptional set, but the corresponding derived category of
\emph{graded} modules is not, because its Grothendieck group has
infinite rank (the set of irreducibles has a free action of $\Zet$
by shifts of grading, hence is infinite).

Let $\D$ be a $\k$-linear triangulated category equipped with a
triangulated auto-equivalence which we call shift of grading (or
just shift), and denote $\F\mapsto \F(1)$; its powers are denoted
$\F\mapsto \F(n)$. For $X,Y\in \D$ set $\Hom_{gr}^i(X,Y)= \oplusl_n
\Hom^i(X,Y(n))$. Assume that $\D$ is of \emph{graded finite type},
i.e., $\Hom^\bu_{gr}(X,Y)$ is finite dimensional for all $X,Y\in
\D$.

An ordered set $\nab^i\in \D$ is {\it graded exceptional} if
$\Hom^\bu_{gr}(\nab^i,\nab^j)=0$ for $i<j$, and
$\Hom^\bu_{gr}(\nab^i,\nab^i)=\k$ for all $i$. We will (slightly
abusing the terminology) say that a graded exceptional set $\nab^i$
is generating if $\D=\bl \{\nab^i(k)\ |\ i\in I, \, k\in \Zet\}
\br$.

 We record the
graded analogs of the previous two Propositions.

The dual set $\del_i$ to a graded exceptional set is defined by
conditions
$$
\Hom^\bu_{gr} (\del_n, \nab^i)=0 \ \ {\rm for}\ \ n>i;$$
$$\del_n \cong \nab^n \mod \D_{<n}^{gr},$$
where $\D_{<n}^{gr}=\bl \nab^i(k)\ |\ i<n, \ k\in\Zet\br$.

Again, it is easy to see  that
$\Hom^\bu_{gr}(\del_i,\nab^i)=\k$, and
$\Hom^\bu_{gr}(\del_i,\nab^j)=0$ for $i\ne j$, that the dual set
equipped with the opposite ordering is exceptional and that it is
unique if it exists.

\begin{Prop}
Let $\nab\subset Ob(\D)$ be a finite graded exceptional set.

a) The triangulated
 subcategory $\C=\bl \nab\br$ generated by $\nab(k)$, $k\in \Zet$
 is both left and right admissible.

b) The dual graded exceptional set exists.
\end{Prop}

\proofpt (b) follows from (a) as in Proposition \ref{exist_dual}. We
prove (a) by induction in the number of elements in $\nab$. If this
number is one, then the functors
 $X\mapsto \oplusl_{n\in \Zet} \Hom^\bu(\nab^1(n),X)\otimes \nab^1(n)$,
$X\mapsto  \oplusl_{n\in \Zet} \Hom^\bu(X,\nab^1(-n))^*\otimes
\nab^1(n) $, from $\D$ to $\bl\{ \nab^1(k)$, $k\in \Zet\} \br$
 are readily seen to be, respectively,
right and left adjoint to the inclusion.
 This provides the base of the induction. The induction step readily follows
from the next Lemma applied to $\C_1=\bl \nab^1(n)\br $,
$\C_2=\bl \nab^2(n), \dots, \nab^n\br$. \epf

\begin{Lem}
Let $\C_1, \ \C_2$ be  full triangulated subcategories of $\D$,
such that $\C_2\subset \C_1^\perp$; define
a full triangulated subcategory $\C\subset \D$ by
 $\C=\bl \C_1\cup \C_2\br$; thus $[\C]=[\C_1]*[\C_2]$. Assume that
both $\C_1$ and
$\C_2$ are left and right admissible in $\D$. Then $\C$ is also left and
right admissible in $\D$.
\end{Lem}

\proofpt See \cite{BK}, Proposition 1.12. \epf

 We now assume that $\D=\bl \nab(k) \br$, $k\in \Zet$ where
 $\nab=(\nab^i)$, $i\in I$ is a graded exceptional set, where the ordered set $I$
 is either  finite
or isomorphic to the ordered set of positive integers. We let
$\del_i$ be the dual graded exceptional set.

The proofs of the following statements are parallel to the proof of
Proposition \ref{propexcep} (see \cite{nilpokon}) and Lemma
\ref{mutalem} respectively.

\begin{Prop}\label{4}
a) There exists a unique $t$-structure $(\D^{\geq 0}, \D^{<0})$ on
$\D$, such that $\nab^i(k)\in \D^{\geq 0}$, $\del^i(k)\in \D^{\leq
0}$ for all $k\in \Zet$, $i\in I$. Moreover, $\D^{\geq 0}$,
$\D^{<0}$ are given by
 $$ \D^{\geq 0}=\< \{\nab^i(k)[d]\ ,\ i\in I, \; d\leq 0,\ k\in\Zet \}
 \>;$$
$$\D^{< 0}=\< \{\del_i(k)[d]\ ,\ i\in I, \; d> 0, \ k\in \Zet\}
\>.$$

b) The $t$-structure is bounded.

c) For $X\in \D$ we have: $X\in  \D^{\geq 0}\iff
\Hom_{gr}^{<0}(\del_i, X)=0$ $\forall i$;

 $X\in D^{<0}\iff
\Hom^{\leq 0}_{gr}(X, \nab^i)=0$ $\forall i$.

d) Let $\A$ denote the heart of the above $t$-structure. Then every
object of $\A$ has finite length. For every $i$, $k$  the image
$L_i(k)$ of the canonical arrow $\tau_{\geq 0}(\del_i(k))\to
\tau_{\leq 0}(\nab^i(k))$ is irreducible. The objects $L_i(k)$ are
pairwise non-isomorphic, and every irreducible object in $\A$ is
isomorphic to $L_i(k)$ for some $i$, $k$. The $t$-structure induces
a $t$-structure on $\D_i=\bl\{ \nab_1(k),\dots, \nab_i(k) \} \br$,
$k\in \Zet$; the heart of the latter, $\A_i$, is the Serre
subcategory in $\A$ generated by $L_1(k),\dots L_i(k)$, $k\in \Zet$.
The map $\tau_{\geq 0}(\del_i(k))\to L_i(k)$ is a projective cover
of $L_i(k)$ in $\A_i$, and $L_i(k)\to \tau_{\leq 0}( \nab^i)$ is an
injective hull of $L_i(k)$ in $\A_i$. \epf

\end{Prop}

\begin{Lem}\label{mutalemgr} Let $(I,\preceq)$
be an ordered set, and $\nab^i\in \D$, $i\in I$ be a graded
exceptional set. Let $\leq$ be another order on $I$; we assume that
either $I$ is finite, or the ordered set $(I,\leq)$ is isomorphic to
$\Zet_{>0}$. Set $\D_{\leq i} =\bl\{ \nab^j(k)\ |\ j\leq i,\, k\in
\Zet\} \br$; $\D_{< i} =\bl \{ \nab^j(k)\ |\ j< i,\,k\in \Zet \}
\br$.

a) For $i\in I$ there exists a unique (up to a unique isomorphism)
object $\nab^{i}_{mut}$ such that $\nab^{i}_{mut}\in \D_{\leq i}\cap
\D_{<i}^\perp$, and $\nab^i_{mut}\cong \nab^i\mod \D_{<i}$.

b) The objects $\nab^i_{mut}$ form a graded exceptional set indexed
by $(I,\leq)$.

c) We have $\D_{\leq i} =\bl \{ \nab^j_{mut}(k)\ |\ j\leq i,\, k\in
\Zet\} \br$, $\bl \nab\br=\bl \nab_{mut} \br$. \epf
\end{Lem}

Thus the notion of a mutation of a graded exceptional set is
defined.

\subsubsection{Tilting objects}\label{tiltobse}
Let $\D$ be a triangulated category  generated by an exceptional
set $\nab$, let $\del$ be the dual set, and let $\A$ be the heart of
the corresponding $t$-structure.

Recall that an object $X\in \A$ is called \emph{tilting} if it has a
filtration with associated graded $N^i=\tau_{\leq 0}(\nab^i)$ and
also has a filtration with associated graded $M_i=\tau_{\geq
0}(\del_i)$.

Classification of indecomposable tilting objects follows from
\cite{tilt} (cf. also \cite{BBM}): for every $i\in I$ there exists a
unique (up to an isomorphism) indecomposable tilting object $T_i\in
\A$ which lies in $\A_{\leq i}$ but not in $\A_{<i}$; every
indecomposable tilting object is isomorphic to $T_i$ for some $i\in
I$. We have a surjective morphism $T_i\to \nab^i$ whose kernel has a
filtration with associated graded $\nab^j$, $j<i$; and an injective
morphism $\del_i\to T_i$ whose cokernel has a filtration with
associated graded $\del_j$, $j<i$.

\begin{Lem}\label{tilt_crite} Assume that
$\Hom^{<0}(\nab^i,\nab^j)=0=\Hom^{<0}(\del_i,\del_j)$.

Let $X\in \D$ be an object, such that
$\Hom^{>0}(\del_k,X)=0=\Hom^{>0}(X,\nab_k)$ for all $k$. Then $X$
lies in $\A$ and it is a tilting object therein.
\end{Lem}
\proofpt The first condition shows that $\del_i, \, \nab^i\in \A$
for all $i$; thus $M_i=\tau_{\geq 0}(\del_i)=\del_i$ and
$N^i=\tau_{\leq 0}(\nab^i)=\nab^i$ .

The second one implies by a standard argument that $X\in \<
\nab_i[n]\> \cap \< \del_i[-n]\> $, $i\in I$, $n\in \Zet_{\geq 0}$.
In particular, $X\in \D^{\leq 0}\cap \D^{\geq 0}=\A$; thus
$\Hom^i(\del_k,X)=0=\Hom^i(X,\nab_k)$ for $i\ne 0$. It is well known
that the latter condition implies that $X\in \langle \del_i\rangle
\cap \langle \nab^i\rangle$; since $M_i=\del_i$, $N^i=\nab^i$ we see
that $X$ is
tilting.  \epf 

\begin{Rem}\label{onpur}
Notice that we have deduced vanishing of $\Hom^i(\del_k,X)$,
$\Hom^i(X,\nab_k)$ for $i<0$ from the corresponding vanishing for
$i>0$. This elementary argument acquires interesting consequences
when it is taken together
 with a derived Koszul equivalence (see  Introduction, discussion
 before the Conjecture in
 \S \ref{koszul_sect}). This is illustrated by the methods of the present
paper. Namely, we will see later that
for a certain class of objects $X\in D^b(\Uqmod)$,
the vanishing of  $\Hom^i(\del_k,X)$,
$\Hom^i(X,\nab_k)$ for $i>0$  follows from {\em purity}
of some irreducible perverse sheaves on the affine flag variety; while
the corresponding vanishing for $i\ne 0$ can be interpreted
as {\em pointwise purity} of these irreducible perverse sheaves.

Recall that purity is a general property of irreducible perverse
$l$-adic sheaves proved in \cite{BBD}, while pointwise purity does
not hold for a general IC sheaf, but is known to hold, e.g., for IC
sheaves of (affine) Schubert varieties. However, Lemma
\ref{tilt_crite} can be used to
 show that
 pointwise purity of irreducible objects
in a category $\A$ of perverse sheaves formally follows from their
purity, provided that a derived Koszul  equivalence between
$D^b(\A)$ and some other (reasonable) derived category is given.
This is an illustration of  the relation\footnote{Pointed out to me
by Victor Ginzburg, cf. \cite{BGS}.} between the pointwise purity
 property and the Koszul duality formalism.

\end{Rem}

\subsection{Sheaves on the nilpotent cone.}\label{CohN}
In this subsection we recall some results and notations from
\cite{nilpokon}.

Till the end of this section \ref{coh} $G$ is assumed to be an
arbitrary semi-simple algebraic group.



For $\lambda\in \Lambda$ we let $\O_{G/B}(\lambda)\in Coh^G(G/B)$
be the corresponding line bundle on $G/B$, and we set
$\O_\lambda=pr^*(\O_{G/B}(\lambda))\in Coh^G(\Nt)$, where
$pr:\Nt\to G/B$ is the projection.

Recall that we consider the action of the multiplicative group $\Gm$
on $\Nt$ given by $t:(\b,x) \mapsto (\b,t^2x)$. Let us equip
$\O_{G/B}(\lambda)$ with the trivial $\Gm$ action, thus we get an
object $\O_{G/B}(\lambda, 0)\in Coh^{G\times \Gm}(G/B)$ (where $\Gm$
acts trivially on $G/B$); and then set
$\O_{\lambda,0}=pr^*(\O_{G/B}(\lambda,0))\in Coh^{G\times
\Gm}(\Nt)$.

The perverse coherent  $t$-structure on $D^G(\N)$ of the middle
perversity\footnote{To simplify notation we use here a
  normalization
different from the one used in \cite{nilpokon}: the two
$t$-structures differ by the shift by $d=\frac{\dim \N}{2}$.}
 is defined by: $\F\in D^{p,\leq 0}(\N)\iff i_O^*(\F)\in
D^{\leq d - d_O}(Coh^G(O))$, $\F\in D^{p,>0}(\N) \iff i_O^!(\F)\in D^{>
d- d_O}(Coh^G(O))$, where $O$ runs over the set of $G$ orbits in
$\N$, and $d_O=\frac{\dim (O)}{2}$, $d=\frac{\dim \N}{2}$.
 All objects in the heart  of this $t$-structure
have finite length. As was pointed out in the Introduction, for a
pair $(O,\LL)$,
 where $O\subset \N$ is a $G$ orbit and $\LL$ is an
irreducible $G$-equivariant vector bundle on $O$, there is a unique
irreducible object $IC_{O,\LL}$ in the heart  supported on the
closure of $O$ and satisfying $IC_{O,\LL}|_{O}=\LL[d_O-d]$.

We will use the bijection (constructed in \cite{nilpokon}) between
$\Lambda^+$ and the set of pairs $(O,\LL)$ as above; we denote it by
$\lambda  \leftrightarrow (O_\lambda, \LL_\lambda)$. For $\lambda\in
\Lambda$ we set $A_\lambda=\pi_*(\O_\lambda)$ (recall that $\pi_*$
denotes the {\em derived} direct image functor).

Set $D^G_{\preceq \lambda}(\N)=\bl \{ A_\nu \ | \ \nu\in
\Lambda^+,\, \nu\preceq \lambda\} \br$, $D^G_{\prec \lambda}=\bl \{
A_\nu \ | \ \nu\in \Lambda^+,\, \nu\prec \lambda\} \br$. The
following facts were proved in \cite{nilpokon}:

\begin{equation}\label{genera}
D^G(\N)=\cupl_\lambda D^G_{\preceq \lambda}(\N);
\end{equation}
\begin{equation}\label{A_perv}
\begin{array}{ll}
 D^{p,> 0}(\N)=\< A_\lambda[n] \ |\ \lambda\in \Lambda^+,\, n< 0\> \\
D^{p,\leq 0}(\N)=\< A_{w_o(\lambda)}[n] \ |\ \lambda\in \Lambda^+,\, n\geq 0\>
\end{array}
\end{equation}
\begin{equation}\label{irrebu}
\dim \Hom(A_{w_o(\lambda)},A_\lambda)=1,\ \ {\rm and }
 \ \ IC_{O_\lambda, \LL_\lambda}\cong Im (A_{w_o(\lambda)}\overset{c_\lambda}{\To} A_\lambda)
\ \ {\rm for}\ \ \lambda \in \Lambda^+,
\end{equation}
where $c_\lambda$ is
 the unique (up to scaling) non-zero morphism, $Im$ denotes the
 image
 with respect to the perverse $t$-structure and $w_o\in W_f$ is the longest
element;
\begin{equation}\label{A_perp}
A_\lambda\in D_{\prec\lambda}^G(\N)^\perp;\ \ \ A_{w_o(\lambda)}\in
^\perp D_{\prec\lambda}^G(\N)\ \ \ \ {\rm for}\ \lambda\in \Lambda^+
\end{equation}
\begin{equation}\label{mod_smaller}
A_\lambda\cong A_{w(\lambda)} \mod D_{\prec\lambda}^G(\N) \ \ \ {\rm
for}\ \ \ \lambda\in \Lambda^+, \ w\in W_f
\end{equation}
\begin{equation}\label{ortohren}
\begin{array}{ll}
D_{\prec\lambda}^G(\N)^\perp = \bl \{ V_\mu \otimes \O_\N\ |\ \mu
\prec
\lambda\} \br ^\perp \\
^\perp(D_{\prec\lambda}^G(\N))= ^\perp \bl \{ V_\mu \otimes \O_\N\
|\ \mu \prec \lambda\} \br.
\end{array}
\end{equation}

Finally, we recall a result of \cite{B}, 
which will be used to deduce Corollary \ref{Humphreys} in section
\ref{sect32} (Humphreys' conjecture) (and will not be used elsewhere
in the paper).

\begin{Prop} Let $\lambda\in \Lambda^+\subset W$, and let $\c\subset W$ be the
two-sided cell containing  $\lambda$. Let $(O_\lambda, \LL_\lambda)$
be the orbit and the irreducible vector bundle on it, corresponding
to $\lambda$ as above. Then $O_\lambda\subset \N$ is the orbit
corresponding to $\c$ under the bijection defined in
\cite{cell4}. \epf
\end{Prop}

\subsection{Exotic $t$-structure.}\label{exotstr}
Some of the arguments appearing in this section can be replaced by
shorter ones relying on results of \cite{AB} or \cite{ABG}. We found
it worthwhile, however, to present direct proofs; though the
intuition comes from \cite{AB}, the techniques employed in this
section go back at least to \cite{Dem}.

The category $ D^{G\times \Gm}(\Nt)$ is equipped with a shift functor
$\F\mapsto \F(1)$, where $\F(1)$ stands for the tensor product of
$\F$ with the tautological character of $\Gm$. For $\F,\G\in
D^{G\times \Gm}(\Nt)$ we have
$\Hom^\bu_{D^G(\Nt)}(\F,\G)=\oplusl_n \Hom^\bu_{D^{G\times
\Gm}(\Nt)} (\F,\G(n)),$ where the restriction of equivariance
functor is omitted from the notation.

 Set $\O_{\lambda, n}=\O_{\lambda, 0}(n)$.

\begin{Lem}\label{no_Hom}
 We have $\Hom^\bu_{D^G(\Nt)}(\O_\lambda,\O_\mu)=0$ unless
$ \lambda\succeq \mu$.
Also, $\Hom^\bu_{D^G(\Nt)}(\O_\lambda,\O_\lambda)=\k$.
\end{Lem}

\proofpt We have to check that $R^\bu\Gamma ^G(\O_\nu)=0$ if
$-\nu\not \in R^+$, and $R^\bu\Gamma ^G(\O)=\k$. Here by $R\Gamma
^G$ we mean the derived functor of the functor $\F\mapsto
\Gamma(\F)^G=\Hom_{Coh^G}(\O,\F)$. We have
$R^i\Gamma^G(\F)=(R^i(\F))^G$, see, e.g., \cite{izvrat}. 

We first prove that for $\F\in D^G(\Nt)$ there is a canonical
isomorphism
\begin{equation}\label{RGH}
R^i\Gamma ^G(\F)\cong (H^i( \n, \F|_\n))^T.
\end{equation}
Recalling that  $\Nt=G\times_B \n$, we see that \eqref{RGH} holds
when $i=0$ and $\F\in Coh^G(\Nt)$. It remains to check that the
$\delta$-functor $\F\mapsto  (H^i( \n, \F|_\n))^T$ is effaceable.
The category $QCoh^G(\Nt)$ of $G$-equivariant quasi-coherent
sheaves has enough objects adjusted to $\Gamma$ (and hence to
$\Gamma^G$), which are of the form $Av(\F)$, $\F\in QCoh(\Nt)$;
here $Av:QCoh(\Nt)\to QCoh^G(\Nt)$ is the "averaging" functor,
right adjoint to the forgetful functor $QCoh^G(\Nt)\to QCoh(\Nt)$.
We have $Av(\F)=a_*pr^*(\F)$, where $pr:G\times \Nt\to \Nt$, and
$a:G\times \Nt\to \Nt$ are the projection and the action maps. It
is easy to see that $\Gamma (Av(\F))|_\n$ is an injective module
over the Lie algebra $\n$. This implies \eqref{RGH}.

We now claim that the right hand side of \eqref{RGH} vanishes for
$\F=\O_\lambda$ if $-\lambda\not \in R^+$. Indeed, the standard
complex for computation of $H^\bu(\n, \Gamma(\O(\lambda))|_\n$ is
$\Lambda(\n^*)\otimes Sym(\n^*)$ with some differential; the action
of $T$ is the natural action twisted by $\lambda$. If $-\lambda\not
\in R^+$, then the space of $T$ invariants in the complex vanishes,
hence so does the space of $T$ invariants in its cohomology.
 \epf

\begin{Prop}\label{equip} Equip $\Lambda$ with any order compatible with
the partial order $\preceq$. Then
the set of objects $\O_\lambda$
 indexed by $\Lambda$ with this order is
an  exceptional set generating $D^G(\Nt)$.

The set  $\O_{\lambda,0}$ is a graded exceptional set in
$D^{G\times \Gm}(\Nt)$, and $D^{G\times \Gm}(\Nt)=\bl \O_{\lambda,
n}\br$, $\lambda\in\Lambda$, $n\in \Zet$.
\end{Prop}

\proofpt The set in question is (graded) exceptional by Lemma
\ref{no_Hom}. The set $\O_\lambda$ generates $D^G(\Nt)$ by
\cite{nilpokon}, Corollary 2 on p.13;
 $\O_{\lambda, n}$ generate $D^{G\times \Gm}(\Nt)$
by a similar argument. \epf

We now introduce another partial
 ordering $\leq$ on $\Lambda$. To this end, recall
the (2-sided) Bruhat partial order on the affine Weyl group $W$.
For $\lambda\in \Lambda$ let $w_\lambda$ be the minimal length
representative of the coset $W_f\lambda\subset W$.


 We set
$\mu\leq \lambda$ if $w_\mu$ precedes $w_\lambda$ in the Bruhat
order.

We will use the following well-known properties of the Bruhat order:

\begin{itemize} \item The two orders $\preceq$, $\leq$ coincide on
$\La^+$.

\item The two orders $\preceq$, $\leq$ coincide on any $W_f$ orbit
in $\La$.

\item If $\la\leq \mu$, then $\la\in conv(\mu)$, where
$conv(\mu)$ denotes the convex hull of the $W_f$ orbit of $\mu$
intersected with $\La$.

\end{itemize}

We fix a complete order
%
 $\leq_{compl}$ on $\Lambda$ compatible
with $\leq$.

 In view of the above
properties of the Bruhat order we
 can (and will) assume that $\leq_{compl}$ satisfies
the following two
requirements.\footnote{One can use the comparison
with geometry of affine flags,
see section \ref{onAB}, to show that the second requirement is, in fact,
redundant. We still found it convenient to impose it for the sake of
a technical simplification.}

\begin{itemize}\item
$(\Lambda, \leq_{compl})$ is isomorphic to $\Zet_{\geq 0}$.

\item Let $conv^0(\mu)$ denote the complement to the $W_f$ orbit
of $\mu$ in $conv(\mu)$. Then we have

\begin{equation}\label{cond2}
\la \in conv^0(\mu)\ \ \Rightarrow \ \ \la \leq _{compl}\mu.
\end{equation}
\end{itemize}


Proposition \ref{equip} implies that the set $\O_\la$ (respectively,
$\O_{\la, 0}$) equipped with any complete order compatible with
$\preceq$  is a (graded) exceptional generating set.
  We define the (graded) exceptional set $\nab_\lambda$
(respectively, $\nab_{\lambda,0}$) to be the $ \leq_{compl}$
mutation of the set $\O_\lambda  $ (respectively, $\O_{\lambda,0}
$). We let $\del_\lambda$ (respectively, $\del_{\lambda, 0}$) be the
dual (graded) exceptional sets. It is clear from the definitions
that the forgetful functor $D^{G\times \Gm}(\Nt)\to D^G(\Nt)$ sends
$\nab^{\lambda,0}$ to $\nab^\lambda$ and $\del_{\lambda,0}$ to
$\del_\la$.

Our next goal is to get a more explicit description of $\nab^\la$,
$\nab^{\la,0}$.
 It is
provided by the next Proposition \ref{indepePr}; to state it we need
some notations.

 For a subset $S\subset \Lambda$ set $D_S^G(\Nt)=\bl \O_\lambda\ |\
\lambda\in S \br$, $D_S^{G\times \Gm}(\Nt)=\bl \O_{\lambda,n}\ |\
\lambda\in S,\, n\in \Zet \br$.

 Let
$\Pi_\la^l:D^G(\Nt)\to D^G_{conv^0(\lambda)}(\Nt)^\perp$,
$(\Pi_\la^l)^{gr}:D^{G\times \Gm}(\Nt)\to
D_{conv^0(\lambda)}^{G\times \Gm}(\Nt)^\perp$ be the functors left
adjoint to the embedding functors.


 For  a simple root $\al$ let $s_\alpha\in W_f$ be the corresponding
simple reflection. Let $\pr_\alpha:G/B\to G/P_{\alpha}$ be the
projection, where $P_\alpha$ is the minimal parabolic corresponding
to $\alpha$.
 Set $\Ntil_\alpha=T^*(G/P_\alpha)\times_{G/P_{\alpha}}
G/B$; the differential of $\pr_\alpha$ provides a closed embedding $
i_\al: \Ntil_\alpha \imbed \Ntil$. We let  $\pit_\al$ denote the
projection $\Ntil_\al \to T^*(G/P_\al)$.

Define a functor $F_\al:D^{G}(\Nt)\to D^{G}(\Nt)$ by
$$F_\al:\F\mapsto \left( i_{\al *} \circ \pit_\al^* \circ
\pit_{\al*}\circ i_\al^!(\F(\rho)) \right) (-\rho);$$ a functor
$F_\al^{gr}:D^{G\times \Gm}(\Nt)\to D^{G\times \Gm}(\Nt)$ is defined
by the same formula.\footnote{\label{foosc}The formula defines a
functor on the equivariant derived category only if $\O_\rho$
carries a $G$ equivariant structure; this is true if $G$ is simply
connected, but not in general. However, the derived category
$D^G(\Nt)$, $D^{G\times \Gm}$ is a direct summand in
$D^{G_{sc}}(\Nt)$, $D^{G_{sc}\times \Gm}(\Nt)$, where $G_{sc}$ is
the universal cover of $G$. It is easy to see that the functors
$F_\al$, $F_\al^{gr}$ defined by the above formula for $G_{sc}$
preserve this direct summand, thus the functors are defined for a
general $G$. This observation applies also to $F_\al'$,
$(F_\al')^{gr}$.}

Notice that each of the functors $F_\al$, $F_\al^{gr}$ is defined as
a composition of six functors which can be grouped into three
adjoint pairs: $(\pit_\al^*, \pit_{\al*})$, $(i_\al^*, i_\al^!)$,
$(T_{-\rho}, T_\rho)$, where we let $T_\la$ denote the functor of
twist by $\O_\la$ (respectively, by $\O_{\la,0}$).

We also define the functors $F'_\al$, $(F'_\al)^{gr}$  by $\F\mapsto
\left( i_{\al
*} \circ \pit_\al^! \circ \pit_{\al*}\circ i_\al^*(\F(-\rho)) \right)
(\rho)$.

 Thus $F_\al$ (respectively, $F_\al^{gr}$) is the
composition of $\pit_{\al*}\circ i_{\al}^! \circ T_{\rho}$ with its
left adjoint. In particular, we get a canonical arrow $F_\al
@>{can_\al}>> id$ and similarly for $F_\al^{gr}$.

 We have a canonical arrow $id @>{can_\al'}>> F_\al'$ and
its graded version.

\begin{Prop}\label{indepePr}
a) We have $\nab^\la\cong \Pi_\la^l(\O_\la)$, $\nab^{\la,n}\cong
(\Pi_\la^l)^{gr}(\O_{\la,n})$. In particular, $\nab^\la$,
$\nab^{\la,n}$ do not depend on the choice of $\leq_{compl}$.

b) Assume that $s_\al(\la)\preceq \la$. Then we have  canonical
distinguished triangles:
$$F_\al^{gr}( \nab^{\la,0}) @>{can_\al}>> \nab^{\la,0} \to
\nab^{s_\al(\la),2},$$
$$ \nab^{\la,0}\to \nab^{s_\al(\la),2}
@>{can_\al'}>>  (F'_\al)^{gr}(\nab^{s_\al(\la),2})$$
and similarly with $F_\al^{gr}$, $(F_\al')^{gr}$, $\nab^{\la,0}$,
$\nab^{s_\al(\la),2}$ replaced by $F_\al$, $F'_\al$, $\nab^\la$,
$\nab^{s_\al(\la)}$.

\end{Prop}

The proof of the Proposition is preceded by two Lemmas.

\begin{Lem}\label{ab}
a) The left adjoint to $F_\al$, $F_\al^{gr}$ is isomorphic to
$F_\al[-2]$, $F_\al^{gr}[-2](-2)$ respectively.

b)
 We have $F_\al\cong F_\al'[2]$, $F_\al^{gr}\cong (F_\al')^{gr}[2](-2).$
\end{Lem}

\proofpt  We  use the well-known isomorphisms of $G\times \Gm$
equivariant line bundles on $\Nt$:
\begin{equation}\label{wki}
\begin{split} 
& \O(-\Nt_\al)\cong \O_{\al, -2},\\
& K_{G/B}\otimes \pi_\al^*(K_{G/P_\al}^{-1})\cong \O_{-\al},
\end{split}\end{equation} where $K$ stands for the canonical line bundle of
top degree forms. They imply canonical isomorphism:
%
$$
 i_\al^!(\F)\cong i_\al^*(\F\otimes
\O(\Nt_\al))[-1]=i_\al^*(\F(-\al, 2))[-1], $$
 $$\pit_\al^!(\F) \cong
\pit^*_\al(\F)\otimes _{\O_{G/B}} \left( K_{G/B}\otimes
\pi_\al^*(K_{G/P_\al}^{-1})\right) [1]\cong
\pit^*_\al(\F) (-\al, 0)[1]. 
$$

Since $\pit_\al^*$, $i_\al^*$ are left adjoint to $\pit_{\al*}$,
$i_\al^*$ respectively, while $\pit_\al^!$, $i_\al^!$ are right
adjoint to $\pit_{\al*}$, $i_\al^*$ respectively, we get (a) by
plugging in the latter isomorphisms. They also imply that
$$F_\al^{gr}\cong T_{-2\rho+\al,0}(F_\al')^{gr} T_{2\rho-\al,2}.$$
Since $\langle \check \al, 2\rho-\al\rangle  =0$ (where $\check \al$
is the dual coroot), we see that $\O_{2\rho-\al}$ is lifted from
$G/P_\al$, hence $T_{2\rho-\al,0}$ commutes with $(F_\al')^{gr}$.
The isomorphism in (b) follows. \epf

\begin{Lem}\label{cd}
a) If $s_\al(\la)\preceq \la$, then we have canonical isomorphisms
in (the quotient category of) $D^{G\times \Gm}(\Nt)$:

\begin{align*}& \O_{s_\al(\la),2} \cong cone
(F_\al(\O_{\la,0})@>{can_\al}>> \O_{\la,0}) \mod
D_{conv^0(\la)}^{G\times \Gm}(\Nt),\\
& \O_{\la,0}[1] \cong cone (\O_{s(\la),2})@>{can'_\al}>>
F'_\al(\O_{s_\al(\la),2}) ) \mod D_{conv^0(\la)}^{G\times \Gm}(\Nt)
\end{align*} and also similar
 isomorphisms in (the quotient category of) $D^{G}(\Nt)$.

 b) The functors $F_\al$, $F^{gr}_\al$ preserve
  $D_{conv(\la)}^G(\Nt)$, $D_{conv(\la)}^{G\times
\Gm}(\Nt)$. 

\end{Lem}

\proofpt We prove (a) and (b) together by induction. We assume that
(b) is known for $\la\in conv^0(\mu)$, and deduce that (a) holds for
$\la \in W(\mu)$. Then the isomorphisms of (a) together with Lemma
\ref{ab}(b) show that (b) also holds for $\la \in W(\mu)$.

If $\langle \check \al ,\la\rangle =0$, then
$i_\al^!(\O_{\la+\rho})=i_\al^*(\O_{\la+\rho-\al})[-1]$ is
isomorphic to the line bundle $\O_{\Pone}(-1)$ when restricted to a
fiber of $\pit_\al$. Thus in this case we have
$F_\al^{gr}(\O_{\la,0})=0$, so we get the desired isomorphism.
Assume now that $\langle \check \al ,\la\rangle >0$.

 For a
vector bundle $\V$ on $\Nt$ we have, in view of \eqref{wki}, a
canonical (in particular, respecting $G$ or $G\times \Gm$
equivariance) short exact sequence
$$0\to \V(\al,-2) \to \V\to i_{\al*} i_\al^* \V \to 0.$$
Applying Serre-Grothendieck duality  (and plugging in $\V^*$ for
$\V$) we get a distinguished triangle
$$i_{\al*}i^!_\al\V \to \V \to \V (-\al,2).$$
If a vector bundle $\cW$ on $\Nt _\al$ is trivial on the fibers of
the projection $\pit_\al$ (if $\cW$ is an equivariant bundle, the
trivialization is not required to be compatible with the equivariant
structure), then we have $\pit_\al^*\pit_{\al*}\cW\iso \cW$. Thus if
$\V$  is a $G\times \Gm$ equivariant vector bundle, such that
$i_\al^!(\V(\rho))[1]=i_\al^*(\V(\rho-\al))$ is trivial on the
fibers of $\pit_\al$, then we have a distinguished triangle $i_{\al
*} \circ \pit_\al^* \circ \pit_{\al*}\circ i_\al^!(\V) \to \V \to
\V(-\al,2)$. If $i_\al^*\V(-\al+\rho)$ is fiberwise trivial, then we
get a distinguished triangle $F_\al(\V) \to \V \to \V(-\al,2)$.

For $\mu\in \La$,  $\langle \check\al ,\mu\rangle \geq 0$, let
$\V_{\mu,n}$ denote the pull-back under the projection $\Nt \to G/B$
of  $\pi_\al^* \pi_{\al*}(\O_{\mu,n})$. It is clear that
$\pi_\al^*\pi_{\al
*}(\O_\mu))$ is trivial on the fibers of $\pi_\al$. It is also easy
to see that this vector bundle carries a filtration with associated
graded $\O_{\mu,n}$, $\O_{\mu-\al,n}, \dots, \O_{s_\al(\mu),n}$.

 Set $\V=\V_{\la-\rho,0}(\rho,0)$, where we have assumed
without lost of generality that $G$ is simply-connected, so
$\O_\rho$ is $G$ equivariant, cf. footnote \ref{foosc}.  Then we
have $i_\al^!(\V)[1]\cong
i_\al^*(\V(-\al))=i_\al^*(\V_{\la-\rho}(-\al+2\rho))$, where we have
dropped the $\Gm$-equivariance to unburden the notation. The vector
bundle  $i_\al^*(\V_{\la-\rho}(-\al+2\rho))$ is trivial on the
fibers of $\pit_\al$, since $\langle \check\al ,-\al+2\rho\rangle
=0$. Thus we have a distinguished triangle
$$F_\al(\V)\to  \V \to \V(-\al,2).$$
On the other hand, the above filtration on $\V_\mu$ yields a
filtration on $\V$ and on $\V(-\al,2)$. In both cases all but one
terms in this filtration lie in $D_{conv^0(\la)}^{G\times
\Gm}(\Nt)$. The only term which does not is $\O_{\la,0}$ for $\V$,
and $\O_{s_\al(\la),2}$ for  $\V(-\al,2)$. Using the induction
assumption we get the first isomorphism for $\la\in W_f(\mu)$.

The second isomorphism is proved similarly, using the observation
that for $\V=\V_{\la-\rho}(\rho-\al)$ and $\cW=i_\al^*(\V)(-\rho)$
we have $\cW\iso \pit_\al^!\pit_{\al*}\cW$.
 \epf

\medskip

{\em Proof of Proposition \ref{indepePr}.} We prove the statements
about $\nab^\la$, the proof of the "graded version" is similar.
  Set
$(\nab')^{\la}=\Pi_\la^l(\O_{\la})$.

We first check that if $s_\al(\la) \prec \la$, then
 we have canonical distinguished triangles
\begin{equation}\label{FF}
F_\al(( \nab')^{\la}) @>{can_\al}>> (\nab')^{\la} \to
(\nab')^{s_\al(\la)},\end{equation}
$$(\nab')^\la \to (\nab')^{s_\al(\la)} @>{can_\al'}>>F_\al'( (\nab')^{s_\al(\la)})
.$$

This follows from Lemmas \ref{ab}, \ref{cd}:   Lemma \ref{cd}(b)
together with Lemma \ref{ab} show that $\Pi_\la^l$ commutes with
$F_\al$, $F_\al'$, then the distinguished triangles follow from
Lemma \ref{cd}(a).

It remains to show that $(\nab')^{\la}\cong \nab^{\la}$.

It suffices to check that
$$
(\nab')^{\la} \cong \O_{\la} \mod D^{G}_{<_{compl} \la},
$$
$$
\Hom^\bu((\nab')^{\la}, (\nab')^{\mu})=0\ \ \ \ \ {\mathrm{if}} \ \
\ \la< _{compl} \mu,$$ where $D^G_{<_{compl} \la}=D_{\{\mu |
\mu<_{compl} \la\} }^G$.

Here the first isomorphism is clear from
the definition of $(\nab')^\la$ and the  requirement \eqref{cond2} on the
order $\leq _{compl}$. The second isomorphism follows from the
following stronger statement:

\begin{equation}\label{notinconv}
\Hom((\nab')^\la,(\nab')^\mu)\ne 0 \ \Rightarrow \ \ \mu \in conv^0
(\la) \bigvee\  (\mu \in W_f(\la)\  \&\  \mu \preceq \la).
\end{equation}

We now prove \eqref{notinconv}. Assume first $\la \in W_f(\mu)$
and  $\mu \not \preceq \la$;
 we show that in this case
\begin{equation}\label{van}
\Hom((\nab')^\la,(\nab')^\mu)=0.
\end{equation}
We have $\Hom ((\nab')^\la,(\nab')^\mu)=\Hom (\O_\la, (\nab')^\mu)$,
because $(\nab')^\la\cong \O_\la\mod D^G_{conv^0(\la)}$ and
$(\nab')^\mu\in  D^G_{conv^0(\la)}(\Nt)^\perp$. On the other hand,
 Lemma
\ref{no_Hom} implies that $(\nab')^\mu\cong \O^\mu \mod D^G_{ \succ
\mu}$,  where $D^G_{ \succ \mu}(\Nt)=D^G_{\{\nu\ |\ \nu \succ \mu\}
}(\Nt)$. Lemma \ref{no_Hom} shows also that $\O_\la\in ^\perp D^G_{
\succeq \mu}(\Nt) $, which yields \eqref{van}.

It remains to show \eqref{van} assuming that $\la \not \in
conv(\mu)$. It is enough to check that $(\nab')^\mu\in
(D^G_{conv(\la)})^\perp$. If $\mu\in \La^+$, then
$(\nab')^\la=\O_\la$, so  $(\nab')^\mu\in (D^G_{conv(\la)})^\perp$
by  Lemma \ref{no_Hom}. Comparing \eqref{FF} with Lemma \ref{cd}(b)
and Lemma \ref{ab}(b) we see that if   $(\nab')^\mu\in
(D^G_{conv(\la)})^\perp$ and $s_\al(\mu )\prec \mu$, then also
$(\nab')^{s_\al(\mu)}\in (D^G_{conv(\la)})^\perp$. Thus
$(\nab')^\mu\in (D^G_{conv(\la)})^\perp$ for all $\mu$ such that
$\la \not \in conv(\mu)$. \epf

\begin{Def}
The exotic $t$-structure on $D^G(\Nt)$ (respectively $D^{G\times
\Gm}(\Nt)$) is the $t$-structure of the (graded) exceptional set
$\nab^\lambda$ (resp. $\nab^{\lambda, 0}$).
\end{Def}

We let $\E^G(\Nt)$, $\E^{G\times \Gm}(\Nt)$ denote the hearts of the
exotic $t$-structures on $D^G(\Nt)$, $D^{G\times \Gm}(\Nt)$
respectively. In view of Propositions \ref{propexcep}, \ref{4} (see
section \ref{qhc}), isomorphism classes of irreducible objects in
$\E^G(\Nt)$ are in bijection with $\Lambda$, while irreducible
objects in $\E^{G\times \Gm}(\Nt)$ are in bijection with
$\Lambda\times \Zet$. For $\lambda\in \Lambda$, $n\in \Zet$ we let
$E_\lambda\in \E^G(\Nt)$, $E_{\lambda,n}\in \E^{G\times \Gm}(\Nt)$
be the corresponding irreducible objects. It is clear from the
construction that the forgetful functor $D^{G\times \Gm}(\Nt)\to
D^G(\Nt)$ sends $E_{\la,n}$ to $E_\la$.

\subsection{Exotic and perverse.}
In this subsection we restrict attention to $D^G(\Nt)$ to simplify
notations; we leave it as an exercise for the reader to provide the
"graded" analogue.

\begin{Lem}\label{nabtoA}
For $\lambda\in \Lambda^+$ and any $w\in W_f$ we have
 $\pi_*(\nab^{w(\lambda)})=A_\lambda$,
$\pi_*(\del_{w(\lambda)})=A_{w_o(\lambda)}$.
\end{Lem}

\proofpt  We prove the first isomorphism, the proof of the second
one is similar.

If $w=1$, then the isomorphism is just the definition of $A_\la$.
Thus we will be done if we show that $\pi_*(\nab^{s_\al(\lambda)})
\cong \pi_*(\nab^\la)$ for any simple reflection $s_\al$ and $\la\in
\La$. This follows from Proposition \ref{indepePr}(b) in view of
$\pi_*\circ F_\al=0$. The latter equality is a consequence of the
fact that for any $\F$ the complex of sheaves $F_\al(\F)$ is
concentrated on $\Nt_\al$ and its restriction to any fiber of
$\pit_\al$ is isomorphic to a sum of several copies of
$\O_{\Pone}(-1)[d]$ for some $d$; thus its derived global sections
vanish. \epf

\begin{Cor}\label{exact}
The functor $\pi_*:D^G(\Nt)\to D^G(\N)$ is $t$-exact, where
$D^G(\Nt)$ is equipped with the exotic, and $D^G(\N)$ with the
perverse coherent $t$-structure of the middle perversity (see
section \ref{CohN}).
\end{Cor}

\proofpt Compare the definition of the exotic $t$-structure with
\eqref{A_perv} and Lemma \ref{nabtoA}. \epf

\begin{Prop}\label{E_IC}
For  $\lambda\in -\Lambda^+$ we have
$\pi_*(E_\lambda)=IC_{O_{w_o(\lambda)}, \LL_{w_o(\lambda)}}$.

For $\la \not \in -\Lambda^+$ we have $\pi_*(E_\la)=0$.
\end{Prop}

\proofpt Recall that for $\la\in \La$ the object $E_\lambda$ is the
image of the unique (up to a constant) non-zero morphism
$C_\lambda:\del_\lambda\to \nab^\lambda$, while for $\la\in \La^+$
the  object $IC_{O_\la,\LL_\la}$ is the image of a unique up to a
constant non-zero morphism $c_\la:A_{w_o(\la)}\to A_\la$ . By Lemma
\ref{nabtoA} and Corollary \ref{exact},
 $\pi_*(E_\lambda)$ is the image of the morphism $\pi_*(C_\lambda):
A_{\lambda}\to A_{w_o(\lambda)}$ if $\lambda\in -\Lambda^+$. Thus it
suffices to show that $\pi_*(C_\la)\ne 0$ iff $\la \in -\La^+$.

Assume first that $\la\in -\La^+$. We have $E_\lambda \cong
\nab_\lambda \mod D_{<_{compl}\la }^G(\Nt)$. Since $\la \in -\La^+$,
we have $\nu<\la \Rightarrow \nu \in conv_0(\la)$; we can also
assume that the complete order $\leq _{compl}$ is chosen so that
$\nu<_{compl}\la \Rightarrow \nu \in conv_0(\la)$. Thus if
$\pi_*(E_\la)=0$, then $\pi_*(\nab^\la)=A_{w_o(\la)}$ lies in the
full triangulated subcategory generated by $\pi_*(\nab^\nu)$,
$\nu\in conv^0(\la)$.  In view of Lemma \ref{nabtoA}, the latter
category is generated by $A_\nu$, $\La^+\owns \nu \prec w_o(\la)$.
This category does not, in fact, contain $A_{w_o(\la)}$ by, e.g.,
\eqref{A_perp}.

Assume now that $\la\not \in -\La^+$, thus there exists a simple
reflection $s_\al$, such that $s_\al(\la)\prec \la$. Consider  (the
non-graded version of) the first distinguished triangle of
Proposition \ref{indepePr}(b). We have $\Hom^\bu(\del_\la,
\nab^{s_\al(\la)})=0$, thus the morphism $C_\la$ factors through a
morphism $\del_\la \to F_\al(\nab^\la)$. It was explained in the
proof of Lemma \ref{nabtoA} that $\pi_*\circ F_\al=0$. Hence in this
case we have $\pi_*(C_\la)=0$. \epf

\subsection{Positivity Lemma.}\label{posisect}
We state the key result, which formally implies our main Theorem
\ref{ET}; the proof appears in section \ref{onAB}.

\begin{Lem}\label{posi_lem}
$\Hom^i(\Delta_{\lambda,n+\delta_\lambda},E_{\mu,\delta_\mu})=0$ if $i
>n$;

$\Hom^i(E_{\mu,\delta_\mu},
\nab^{\lambda,n+\delta_\lambda})=0$
if $i>-n$ (see \ref{nota} for
notations).

\end{Lem}

\begin{Rem}
The Lemma is equivalent to the statement that some $\Zet$-graded
vector spaces have trivial components of negative degree, hence the
name.
\end{Rem}

\begin{Rem}\label{purpur}
There are several results in representation theory, the
 Kazhdan-Lusztig conjecture being the first and the most famous one,
which are proved by identifying an algebraically defined vector
space with a (co)stalk of an irreducible perverse sheaf, and then
using the deep information about the action of Frobenius on this
(co)stalk provided by Purity Theorem from \cite{BBD}. Our strategy
for proving the key Lemma \ref{posi_lem} also follows this pattern.
More precisely, the proof of Lemma \ref{posi_lem} relies on purity
of perverse sheaves on the affine flag variety; see also Remark
\ref{onpur} in section \ref{tiltobse} above.
\end{Rem}

\section{Quantum group modules: proof modulo the Positivity Lemma.}
\label{qua} In this section the semi-simple group $G$ is assumed to
be of adjoint type.

We start by fixing notations. Recall that irreducible objects of
$\Uq{\textrm -mod}$ are parameterized by their highest weight
$\lambda\in \Lambda^+$. For $\lambda\in \Lambda^+$ let $V(\lambda)$
be the corresponding irreducible representation. We have the dot
action of the affine Weyl group $W$ on $\Lambda$, $w:\mu\mapsto
w\cdot \mu =w(\mu+\rho)-\rho$ (in particular, the subgroup
$\Lambda\subset W$ acts by $\lambda:\mu\mapsto \mu+l\lambda$). By
the linkage principle \cite{APW} the irreducible object $V(\lambda)$
lies in the block $\Uqmod$ iff $\lambda$ lies in the $W$ orbit of 0.
For $w\in W$ we have  $w\cdot 0\in \Lambda^+$ iff $\ell(w)\leq \ell
(w_fw)$ for all $w_f\in W_f$. There is exactly one such element in
each left coset of $W_f$ in $W$. On the other hand, $\Lambda\iso
W_f\backslash W$. Recall that $w_\lambda$ is the minimal length
representative of the coset $W_f\lambda\subset W$, and set
$L_\lambda = V(w_\lambda\cdot 0)$. For example, if $\lambda\in
\Lambda ^+$ we have $L_\lambda= V(l\lambda)$, while for a strictly
anti-dominant $\lambda$ we have $L_\lambda=V(w_o(\lambda)-2\rho)$,
where $w_o\in W_f$ is the longest element.

Similarly, we let $\Weyl_\lambda=\Weyl(w_\lambda\cdot 0)$,
$\coWeyl^\lambda=\coWeyl(w_\lambda\cdot 0)$ be the Weyl module
and dual Weyl module with highest weight $w_\lambda\cdot 0$, and
$T_\lambda=T(w_\lambda\cdot 0)$
be the corresponding tilting module.

\subsection{Some results of \cite{ABG}.}\label{ABGsect}
The result of \cite{ABG} yields a functor $\Psi:D^{G\times
\Gm}(\Nt) \to D^b(\Uqmod)$.

It satisfies the following properties
\begin{equation}\label{twist}
\Psi(\F(1))\cong \Psi(\F)[ 1] 
\end{equation}

\begin{equation}\label{equiv}
\Psi:\oplusl_{n\in \Zet} \Hom(\F, \G[n](-n)) \iso \Hom (\Psi(\F), \Psi(\G));
\end{equation}

\begin{equation}\label{Olala}
\Psi:\O_{\lambda,0}\mapsto RInd_\Bq^\Uq(l\lambda);
\end{equation}

\begin{equation}\label{b_cohomo}
H(b_q,\Psi(\F))\cong \Gamma (\F|_\n),
\end{equation}
the latter isomorphism is $B$-equivariant (here we use notation
$\Gamma$ for the functor of global sections on $D^{B\times \Gm}(\n)$ and
 we write ``$\Gamma$'' rather than ``$R\Gamma$'' since $\n$ is affine,
thus the functor of global sections  is exact).
Also, for $V\in Rep(G)$ we have
\begin{equation}\label{Fr_twist}
\Psi(\F\otimes V)\cong \Psi(\F) \otimes Fr^*(V),
\end{equation}
where $Fr:\Uq\to U(\g)$ is Lusztig's quantum Frobenius morphism.

\medskip

Recall the functors $H,\Ht$ from the Introduction.

\begin{Prop}\label{u_cohomo}
$H(u_q,\Psi(\F))\cong R\Gamma(\F)=\Gamma (\pi_*(\F)).$
\end{Prop}

\proofpt The left hand side can be rewritten as $\oplusl_{\lambda\in
\Lambda^+} V_\lambda \otimes \Hom^\bu_{\Uq}(\k,\Psi(\F)\otimes
Fr^*(V_\lambda^*))$ where $V_\lambda$ is an irreducible $U(\g)$
module with highest weight $\lambda$ (see \cite{AG}). Using
\eqref{Olala} for $\lambda=0$, and \eqref{Fr_twist} we can rewrite
it as $\oplusl_\lambda  V_\lambda \otimes
\Hom^\bu_{D^G(\Nt)}(\O,\F\otimes V_\lambda^*)= \Hom^\bu_{D(\Nt)}(\O
,\F)=R\Gamma(\F)$. \epf
 
We are ready to give the precise statement of the main result.

\begin{Thm}\label{ET}
$\Psi(E_{\lambda,\delta_\lambda})\cong T_\lambda$.
\end{Thm}

\subsection{Consequences of Theorem \ref{ET}.}\label{sect32}
 Recall from section \ref{CohN} perverse IC
sheaves $IC_{O,\LL}\in D^G(\N)$ (introduced in \cite{nilpokon}) and
the bijection $\lambda \leftrightarrow (O_\lambda, \LL_\lambda)$
between $\Lambda^+$ and pairs $(O,\LL)$ where $O\subset \N$ is a $G$
orbit and $\LL$ is an irreducible $G$-equivariant vector bundle
(introduced in \cite{cell4}, see also
 \cite{nilpokon}, \cite{B}).

The following result contains  a more precise version of Theorem
\ref{cohoti}.

\begin{Cor}\label{HTlambda}
a) For $\lambda\in \Lambda$ the sheaf $\Ht(T_\lambda)$ is isomorphic to the
direct sum of all cohomology sheaves of the complex $E_\lambda$.

b) For $\lambda \in -\Lambda^+$ we
 have $$H(T_\lambda)\cong R\Gamma (E_\lambda)\cong \Gamma (IC_{O_{w_o(\lambda)}},
\LL_{w_o(\lambda)}).$$ For $\lambda\not \in - \Lambda^+$ we have
 $$H(T_\lambda)\cong R\Gamma (E_\lambda)=0.$$

\end{Cor}
\proofpt a) follows directly from Theorem \ref{ET} and
\eqref{b_cohomo}. The isomorphism $H(T_\lambda)\cong R\Gamma
(E_\lambda)$ follows from Theorem \ref{ET} and Proposition
\ref{u_cohomo}. Then (b) follows from Proposition \ref{E_IC}. \epf

\begin{Rem}
The fact that $H(T_\lambda)=0$ for $\la\not \in -\La^+$, which was
deduced above from properties of exotic sheaves, is well-known, see,
e.g.,  \cite{O}, \S 3.2.
\end{Rem}

Comparing this statement
with  Proposition 2 of \cite{B} we get the following version
of a conjecture
by J.~Humphreys~\cite{Hu}.

\begin{Cor}\label{Humphreys}
 Let $\lambda\in - \Lambda^+$ be an anti-dominant weight;
let $\c$ be the 2-sided cell in the affine Weyl group $W$ containing
$\lambda$ (where we use the standard embedding $\Lambda \subset W$).
The support of $H(T_\lambda)$ as a coherent sheaf on $\N$ is the
closure of the nilpotent orbit $O_\c\subset \N$, which corresponds
to $\c$ via the bijection of \cite{cell4}. \epf
\end{Cor}

\begin{Rem}\label{proO}
Corollaries \ref{HTlambda}(b) and \ref{Humphreys} imply conjectures in \S 3.2
of
\cite{O}. Indeed, the classes of $IC_{O,\LL}$ in the Grothendieck group
$K(Coh^{G\times \Gm}(\N))$ form the canonical basis of the latter group in the
 sense of \cite{O}; this is proved in \cite{B}.
\end{Rem}

\begin{Ex}\label{Duflo} (cf. \cite{O}, end of \S 3.2)
 When the equivariant vector bundle $\LL$ is trivial
we have $IC_{O,\LL}=j_*(\O_O)[d_O-d] =
Nm_*(\O_{\hat O})[d_O-d]$ where $j:O\imbed \N$ is the embedding, $j_*$
stands for the \emph{non-derived} direct image, and $Nm: \hat O\to \bar O$
 is the
normalization morphism for  the closure $\bar O$ of $O$
(see \cite{nilpokon}, Remark 10). Thus for some indecomposable tilting
$T_\lambda$ we get $H(T_\lambda)=\O_{\hat O}$.
 In fact,  \cite{B}, Proposition 1
implies that this happens precisely when $\la\in -\Lambda^+$ is such
that the minimal length representative of the two-sided coset 
$W_f\la W_f\subset W$ is a Duflo involution.

\end{Ex}

\subsection{Proof of Theorem \ref{ET}  modulo the Positivity Lemma.}
\label{33}
We start with

\begin{Lem}\label{Psiotnab}
 $\Psi(\nab_{\lambda,0}) \cong \coWeyl ^\lambda [-\delta(\lambda)]$,
 $\Psi(\del_{\lambda,0})\cong \Weyl_\lambda[-\delta(\lambda)]$.
\end{Lem}

\proofpt We prove the first isomorphism, the second one is similar.

Set $D^b_{< \lambda}(\Uqmod)=\bl L_\mu \ |\ \mu < \lambda
\br$, and define
 $D^b_{<_{compl} \lambda}(\Uqmod)$  similarly.  
 Recall the Quantum Weak Borel-Weil Theorem of \cite{APW},
which implies that
$$RInd_\Bq^\Uq(l\lambda)\cong
L_\lambda[-\delta_\lambda] \cong \coWeyl^\lambda [-\delta_\la] 
 \mod D^b_{<
\lambda}(\Uqmod).$$
 Using \eqref{Olala} we see that
%
$\Psi(\O_{\lambda,0})\cong \coWeyl^\lambda[-\delta(\lambda)] \mod
D^b_{< \lambda}(\Uqmod).$

It follows by induction in $\la$ that $D^b_{\leq \lambda}(\Uqmod)=\Psi
(D^{G\times \Gm}_{\leq \la}(\Nt))$, where $D^{G\times \Gm}_{\leq \la}(\Nt)
=\bl \O_{\mu,n}\ |\ \mu < \la \br$, and similarly with $\leq$ replaced
 by $\leq_{compl}$ or $< _{compl}$.
 Recall that by the definition of $\nab^{\la,n}$ we have
$$\nab^{\la,0}\in D_{<_{compl}} ^{G\times \Gm}(\Nt)^\perp,$$
$$\nab^{\la,0}\cong \O_{\la,0} \mod D_{<_{compl}} ^{G\times \Gm}(\Nt).$$
In view of 
\eqref{equiv}
this implies: $$\Psi(\nab^{\la,0})\in D_{<_{compl}} (\Uqmod)^\perp,$$
$$\Psi(\nab^{\la,0})\cong L_\la[-\delta_\la]\cong \coWeyl^\la[-\delta_\la]
\mod  D_{<_{compl}}^b(\Uqmod) .$$
It follows that both $\Psi(\nab^{\la,0})[\delta_\la]$ and
$\coWeyl^\la$ are injective hulls of $L_\la$ in the Serre subcategory
of $\Uqmod$ generated by $L_\mu$ with $\mu \leq_{compl}\la$; hence
these two objects are isomorphic.
 \epf

We are now ready to prove Theorem \ref{ET}.

Comparing the Positivity Lemma \ref{posi_lem} with \eqref{twist}
and \eqref{equiv} we see that 
$$\Hom^k(\Weyl_\mu, \Psi(E_{\lambda,\delta_\lambda}))=\oplusl_{i+n=-k}
\Hom^i(\del_{\mu,n+\delta_\mu}, E_{\lambda,\delta_\lambda})=0$$ for $k>0$, and
similarly $\Hom^k(\Psi(E_{\lambda,\delta_\lambda}), \coWeyl^\mu)=0$
for $k>0$. By
Lemma \ref{tilt_crite}
this implies that $\Psi(E_{\lambda,\delta_\lambda})$ is a
tilting object of $\Uqmod$.

The endomorphism algebra $\End (\Psi(E_{\lambda,\delta_\lambda}))
\cong \Hom^\bu( E_\lambda,E_\lambda)$ carries a grading by
non-negative integers, with a one-dimensional component of degree
zero (namely, the homological grading). Hence it has no nontrivial
idempotents, thus $\Psi(E_{\lambda,\delta_\lambda})$ is
indecomposable. Finally, (the proof of) Lemma \ref{Psiotnab} implies that
 $\Psi(E_{\lambda,\delta_\lambda})\in D^b_{\leq_{compl} \lambda}(\Uqmod)$,
but  $\Psi(E_{\lambda,\delta_\lambda})\not
 \in D^b_{<_{compl} \lambda}(\Uqmod)$, hence
$\Psi(E_{\lambda,\delta_\lambda})\cong T_\lambda$.  \epf

\section{Constructible sheaves on affine flags:
proof of Positivity Lemma.}\label{onAB} In this section $G$ is an
arbitrary semi-simple group.

The objective of this section is to express the Hom spaces in
Lemma \ref{posi_lem} as eigenspaces of  Frobenius acting on a
(co)stalk of an IC sheaf; this is achieved in Corollary
\ref{graded_stalk}.

\subsection{Perverse sheaves on affine flags for $\LG$.}
 Recall some  results  of \cite{AB}.\footnote{Notice
a discrepancy in notations between \cite{AB} and the present
paper: the group denoted by $G$ here corresponds to $\LG$ of
\cite{AB}, and vice versa.} Fix $q=p^a$ where $p$ is a prime
number. Consider the Langlands dual group $\LG$ over the base
field $\Fqbar$; the ``loop'' group ind-scheme $\LGK$ and its group
subschemes $\LGO\supset \bI, \, \bI^-$; thus we have
$\LGK(\Fqbar)=\LG\left( \Fqbar((t))\right)$, $\LGO(\Fqbar)=\LG(O)$,
$\bI(\Fqbar)=I$, $\bI^-(\Fqbar)=I^-$ where $O=\Fqbar[[t]]$, and
$I$, $I^-$ are  opposite Iwahori subgroups in $\LG\left(
\Fqbar((t))\right)$. The
affine flag variety (ind-scheme) $\Fl$ is the homogeneous space
$\LGK/\bI$.

\subsubsection{Iwahori-Whittaker sheaves and equivalence $\Phi$}
\label{IWP}
Let  $\P_{asph}$ denote the
 \emph{category of Iwahori-Whittaker perverse sheaves}; this
is a  subcategory in the category
of perverse $l$-adic sheaves on $\Fl$ defined in \cite{AB},
  \S 2.5 (it is denoted $\PIW$ in \emph{loc. cit.}); it consists of
perverse sheaves, which ``transform as a non-degenerate character''
under the action of the pro-unipotent radical  of the Iwahori group
$\bI^-$ . The group scheme $\LGO$
 acts on $\Fl$, and the orbits of
this action are indexed by the set $W_f\backslash W\isol \Lambda$.
For $\lambda\in \Lambda$ let $\Fl^\lambda$ be the corresponding
$\LGO$ orbit, and $i_\lambda:\Fl^\lambda\imbed \Fl$ be the
embedding.

For $\lambda\in \Lambda$ there is a unique (up to an isomorphism)
irreducible object $IC_\lambda^{et}\in \P_{asph}$ whose support is
the closure of $\Fl^\lambda$. Also there is a unique standard object
$\del_\lambda^{et}$ (denoted by  $\del_{w_\lambda}$ in \cite{AB})
characterized by $i_\mu^*(\del_\lambda^{et})=0$ for $\mu\ne\lambda$,
$i_\lambda^*(\del_\lambda^{et})\cong i_\lambda^*(IC_\lambda^{et})$
and a unique costandard object $\nab^\lambda_{et}$ (denoted by
$\nab_{w_\lambda}$ in \emph{loc. cit.}) characterized by $
i_\mu^!(\nab^\lambda_{et})=0$ for $\mu\ne\lambda$,
$i_\lambda^*(\nab^\lambda_{et})\cong i_\lambda^* (IC_\lambda^{et})$.
(The sub or superscript ``et'' (for etale) is used to remind that we
are in a context different from that of section \ref{coh}). For
$\F\in D^b(\P_{asph})$, and a point $x_\lambda\in \Fl^\lambda$,
which lies in the open orbit of an opposite Iwahori subgroup $\bI^-$
we have
\begin{equation}\label{stalk_costalk}
\begin{array}{ll}
\Hom^\bu(\del_\lambda^{et},\F)= \iota_\lambda^!(\F)[d_\lambda] \\
\Hom^\bu(\F,\nab_\lambda^{et})= (\iota_\lambda^*(\F)[-d_\lambda])^*,
\end{array}
\end{equation}
where $\iota_\lambda:x_\lambda\imbed \Fl$ is the embedding, and $d_\lambda=
\dim \Fl^\lambda$.

The objects $\Delta_\lambda^{et}$, $\nab^\lambda_{et}$ can be
described as extension by zero (respectively, direct image) of an
 Artin-Schreier sheaf on $\Fl^\lambda$; for future reference we spell
this out for $\lambda=0$.
 Let $\LB$, $\LB_-$ be the
images  under the projection $\LGO\to \LG$
of the Iwahori subgroups $\bI$, $\bI^-$ respectively; let
$\LN_-\subset \LB_-$ be the unipotent radical. Let $j:\LN_-\imbed \LG/\LB$
be given by $n\mapsto n\LB$; thus $j$ is an open embedding.
Let $\psi:\LN_-\to \Aone$ be a non-degenerate character.
Define the Whittaker sheaf on $\LG/\LB$ by
$$Wh=j_*\psi^*(AS)[\dim \LG/\LB]\cong j_!\psi^*(AS)[\dim \LG/\LB],$$
where $AS$
denotes the Artin-Schreier sheaf. Then
 $\Delta_0^{et}\cong \nab^0_{et}$ is the direct image of $Wh$ under the closed
embedding $i_0:\Fl^0=\LG/\LB\imbed \Fl$.

\medskip

One of the main results of \cite{AB} is  an equivalence of
triangulated categories
$$\Phi:D^G(\Nt)\cong D^b(\P_{asph}). $$

\subsubsection{$\bI$-equivariant sheaves} Let $D_\bI(\Fl)$ be the
 $\bI$-equivariant derived category of $l$-adic sheaves on $\Fl$,
and $\P_\bI\subset D_\bI(\Fl)$ be the subcategory of perverse
sheaves. Convolution of equivariant complexes provides
$D_\bI(\Fl)$ with a structure of monoidal category, and it
provides $D^b(\P_{asph})$ with a structure of a (right) module
category;\footnote{Here we use that $D^b(\P_{asph})$ is identified
with a full subcategory in the derived category of constructible
sheaves on $\Fl$, cf. \cite{AB}, Lemma 1.}
 we denote convolution of
complexes $\F\in D^b(\P_{asph})$, $\G\in D_\bI(\Fl)$ by $\F*\G$. 

The equivalence $\Phi$ is obtained as a composition
$\Phi=Av_{\bI^-,\psi}\circ F$. Here $F$ is a certain
(monoidal) functor $D^G(\Nt)\to D_\bI(\Fl)$
and $Av_{\bI^-,\psi}$ is the
functor of "averaging against a character" $\psi$, i.e.,
$Av_{\bI^-,\psi}:\F\mapsto \del_0^{et}*\F$.

We have
$F(\O_\lambda)=J_\lambda$ where $J_\lambda$ is the \emph{Wakimoto
sheaf}, see \cite{AB},  \cite{ABG}, \cite{Haines}
and $F(V_\lambda\otimes \O_{\Nt})=Z_\lambda$ is the {\em central
sheaf} constructed in \cite{zentr}.

 Also, it
is proved in \cite{AB} that $Av_{\bI^-,\psi}:\P_I\to \P_{asph}$,
i.e. $Av_{\bI^-,\psi}$ is $t$-exact.

The $\bI$ orbits in $\Fl$, also called Schubert cells, are indexed
by $W$; for $w\in W$ let $\Sch_w$ denote the corresponding
Schubert cell, and $j_w:\Sch_w\imbed \Fl$ be the embedding. We
have $J_\lambda\in \P_I$, the support of $J_\lambda$ is the
closure of $\Sch_\lambda$, and we have
\begin{equation}\label{AvWakimot}
\begin{array}{ll}
supp (Av_{I^-,\psi}(J_\lambda))=\overline{\Fl^\lambda}
\\
Av_{I^-,\psi}(J_\lambda)|_{\Fl^\lambda}\cong
\del_\lambda^{et}|_{\Fl^\lambda}, 
\end{array}
\end{equation}
where $supp$ stands for "support", and $\overline{\Fl^\lambda}$
is the closure of $\Fl^\lambda$.
 For $w\in W$ set $j_{w!}=j_{w!}(\cons[\ell(w)])$,
$j_{w*}=j_{w*}(\cons[\ell(w)])$.

We will need another simple geometric property of the Wakimoto
sheaves.

\begin{Lem}\label{Waki_shriek}
We have $i_\lambda^*(J_\lambda)\cong i_\lambda^*(j_{\lambda!}).
$
\end{Lem}

\proofpt It suffices to show that $\Hom^\bu(J_\lambda, j_{w*})=0$
for $w\in W_f\lambda$, $w\ne \lambda$. Let $w_{max}\in W$ be the
maximal length element of the coset $W_f\lambda$, and set
$w_1=w_{max}w^{-1}$. We have $w\in W_f$ and
$\ell(w_1w)=\ell(w_1)+\ell(w)$. It follows that
$j_{w_1}*j_{w*}=j_{w_{max}*}$. Since the functor $\F\mapsto
j_{w_1*}*\F$ is an equivalence (with inverse equivalence given by
$\F\mapsto j_{w_1!}*\F$), it is enough to check that
$\Hom^\bu(j_{w_1*}*J_\lambda, j_{w_{max}*})=0$. Thus we will be done
if we check that the (closure of the) support of $j_{w_1*}*J_\lambda$
does not contain $\Sch_{w_{max}}$. We claim that in fact this
support coincides with the closure of $\Sch_{w_1\lambda}$, which
clearly does not contain $\Sch_{w_{max}}$. To check that
$supp(j_{w_1*}*J_\lambda)=\overline{\Sch_{w_1\lambda}}$ we check
that $j_{w_1*}*J_\lambda$ is a perverse sheaf, and that the Euler
characteristic of its stalk at a point $x\in \Fl$ equals 1 if
$x\in \Sch_{w_1\lambda}$, and 0 otherwise. The claim about the
Euler characteristic is standard. 
 To prove that
$j_{w_1}*J_\lambda$ is perverse we use a trick due to Mirkovi\' c
(cf. \cite{AB}, Appendix 6.1).
 Recall that if $\lambda=\mu_1-\mu_2$, where
$\mu_1,\mu_2\in \Lambda^+$, then
$J_\lambda=j_{\mu_1*}*j_{-\mu_2!}$. Since $w_1\in W_f$,
$\mu\in\Lambda^+$ we have $\ell(w_1\mu_1)=\ell(w_1)+\ell(\mu_1)$.
Thus we have
$j_{w_1*}*J_\lambda=j_{w_1*}*j_{\mu_1*}*j_{-\mu_2!}=j_{w_1\mu_1*}*j_{-\mu_2!}$.
The convolution $j_{u*}*j_{v!}$ is perverse for any $u,v\in W$, see
\cite{AB}, Appendix 6.1. \epf

\subsubsection{Exotic coherent and perverse constructible sheaves}
\label{ecohper}

\begin{Lem}\label{Phinab}
$\Phi(\nab^\lambda)\cong \nab^\lambda_{et}$, $\Phi(\del_\lambda)\cong
\del_\lambda^{et}$.
\end{Lem}

\proofpt Let $D_{\preceq \lambda}^b(\P_{asph})$ (respectively,
$D_{\prec \lambda}^b(\P_{asph})$) be the subcategory of complexes
supported on  $\overline{\Fl^\lambda}$ (respectively, on
$\overline{\Fl^\lambda}\setminus \Fl^\lambda$). In view of
\eqref{AvWakimot} we get, by induction in $\lambda$, that
$$D_{\preceq \lambda}^b(\P_{asph})=\bl
Av_{I^-,\psi}(J_\lambda)\br;$$ since
$Av_{I^-,\psi}(J_\lambda)=\Phi(\O_\lambda)$, we have $D_{\preceq
\lambda}^b(\P_{asph})=\Phi(D^G_{\preceq \lambda}(\Nt))$; and
similarly with $\preceq$ replaced by $\prec$. Thus
$\Phi(\nab_\lambda)\in D_{\prec \lambda}^b(\P_{asph})^\perp\cap
D_{\preceq \lambda}^b(\P_{asph})$, which means that
$\Phi(\nab^\lambda)$ is supported on $\overline{\Fl^\lambda}$ and
its  ``shriek'' restriction to the boundary of $\Fl^\lambda$ vanishes.
Also we have
$$\Phi(\nab_\lambda)\cong \Phi(\O_\lambda) \mod D_{\prec
\lambda}^b(\P_{asph}),$$ which means that
$\Phi(\nab^\lambda)|_{\Fl^\lambda}\cong
Av_{I^-,\psi}(J_\lambda)|_{\Fl^\lambda}$. Thus we are done by
\eqref{AvWakimot}. \epf

\begin{Cor}\label{PhiIC}
$\Phi$ sends the exotic $t$-structure to the tautological one and
$\Phi(E_\lambda)\cong IC_\lambda^{et}$.
\end{Cor}

\proofpt Compare the definition of the exotic $t$-structure with
Remark \ref{izvrat_puch}. Also observe that $E_\lambda$ is the image of
the unique (up to a constant) non-zero morphism $\del_\lambda\to
\nab^\lambda$, while $IC_\lambda^{et}$ is the image of the unique
morphism $\del_\lambda^{et}\to \nab^\lambda_{et}$. \epf

\begin{Cor}\label{Corstalk} For $\F\in D^G(\Nt)$ we have:
$\Hom^\bu(\Delta_{\lambda},\F)\cong \iota^*_\lambda(\Phi(\F))$;
$\Hom^\bu(\F, \nab^{\lambda})^*\cong \iota_\lambda^!(\Phi(\F))$.
\epf \end{Cor}

\subsection{Frobenius weights.}
By a Weil complex  on a scheme $X$ defined over $\Fq$ we will mean
an object  $\F$ in the
``derived category'' of $l$-adic sheaves on $X_{\Fqbar}$
 equipped with an isomorphism
$Fr^*(\F)\cong \F$; a Weil perverse sheaf is a Weil complex, which
is also a perverse sheaf.
 By \cite{BBD}, Proposition 5.1.2, the category of Weil perverse
sheaves on $X$ contains
 the category
of perverse sheaves on the scheme over $\Fq$ as a full subcategory.

 Recall the compatibility
of $F$  with Frobenius. Let $\bq:\Nt\to \Nt$ be the
multiplication by $q$ map, $\bq:(\b,x)\mapsto (\b,qx)$. Then we
have (see [AB], Proposition 1):
\begin{equation}\label{compafr}F\circ Fr^*\cong \bq^* \circ F.
\end{equation}

Fix a square root of $q$, $q^{1/2}\in \Ql$.
Then for $\F\in D^{G\times \Gm}(\Nt)$ the $\Gm$-equivariant structure
induces an isomorphism of objects in $D^G(\Nt)$:
$$\bq^*(\F)\cong \F$$
(recall that
the action of $\Gm$ on $\Nt$ is given by $t:(\b ,x)\mapsto (\b ,t^2x)$).

By means of \eqref{compafr} we get an isomorphism
\begin{equation}\label{FrF}
Fr^*(F(\F))\cong F(\F).\end{equation}

Thus we get a functor $\Ft$ from $D^{G\times \Gm}(\Nt)$ to
the category of Weil complexes on $\Fl$.

Notice that the choice of $q^{1/2}$ defines also a square root
$\F\mapsto \F(\frac{1}{2})$ of the functor of Tate twist on the category
of Weil sheaves (complexes) on an $\Fq$-scheme.

We upgrade the sheaf $\Delta_0^{et}=\nab^0_{et}=i_{0*}(Wh)$
(cf. section \ref{IWP})
 to a Weil sheaf $ \widetilde \Delta_0^{et}= \widetilde \nab^0_{et}=
i_{0*}(\widetilde{Wh})$ where
$$\widetilde{Wh}=j_*\psi^*(AS)[\dim \LG/\LB](\frac{\dim \LG/\LB}{2})\cong j_!\psi^*(AS)[\dim \LG/\LB](\frac{\dim \LG/\LB}{2}).$$

It is clear that $\widetilde\Delta_0^{et}$ is  pure of weight zero.

Then we get a functorial isomorphism $Fr^*(\Phi(\F))\cong \Phi(\F)$
for all $\F\in D^{G\times \Gm}(\Nt)$.
 In particular, for $\F\in
\E^{G\times \Gm}(\Nt)$ the perverse sheaf $\Phi(\F)$ is equipped with a
Weil structure,
 thus we get a functor $\Phit$
from $\E^{G\times \Gm}(\Nt)$ to Weil perverse sheaves on $\Fl$.

 Let us record the following statement, which is immediate from the
 definitions.

\begin{Lem}\label{gr_Fr}
a) For $\F\in \E^{G\times \Gm}(\Nt)$ we have  canonical isomorphisms
$\Ft(\F(n))\cong \Ft(\F)(-\frac{n}{2});$
$\Phit (\F(n)) \cong \Phit(\F)(-\frac{n}{2})$ (see beginning
of section \ref{exotstr} for notations).

b) For $\F,\G \in D^{G\times \Gm}(\Nt)$ the Frobenius action on
$$\Hom^\bu(\Phi(\F), \Phi(\G))=\Hom^\bu_{D^G(\Nt)}(\F,\G) = \oplusl _n
\Hom^\bu_{D^{G\times \Gm}(\Nt)} (\F,\G(n))$$ preserves the direct sum
decomposition, and equals $q^{-n/2}$ on the $n$-th summand. \epf

\end{Lem}

\begin{Lem}\label{Olachist}
Set $J_{\lambda,n}=\Ft(\O_{\lambda,n})$.

Then the Weil perverse sheaf $J_{\lambda,0}$ is pure of weight zero  on
the open dense stratum in its support.
\end{Lem}

\proofpt Let $V_\lambda$ be an irreducible representation of $G$ with
 extremal
weight $\lambda$.
 Recall from \cite{AB} that
 $Z_\lambda=F(V_\lambda\otimes \O_\Nt)$ is
 the
  central sheaf introduced in \cite{zentr}. In particular
it is a perverse sheaf. The sheaf  $V_\lambda\otimes \O_\Nt$
carries an obvious $\Gm$ equivariant structure;
 let $\tilde Z_\lambda=\Ft(V_\lambda\otimes \O_\Nt)$ be the corresponding Weil
perverse sheaf. It follows from the construction of \eqref{FrF}
in \cite{AB} that this Weil structure on $Z_\lambda$
 coincides with the  Weil structure
defined in \cite{zentr}.

The $\Gm$ equivariant vector bundle $V_\lambda \otimes \O_\Nt$
carries a filtration with associated graded $\O_{\nu, 0}$, where
$\nu$ runs over weights of $V_\lambda$. It induces a filtration on
$\tilde Z_\lambda$ with subquotient $J_{\nu,0}$.
 Furthermore, it
is known (see \cite{AB}, \cite{Haines})
 that $\Sch_{\lambda}$ is open in the support of both
$J_\lambda$ and $Z_\lambda$, dense in the support of $J_\lambda$
and does not intersect
the support of $J_\nu$, if $\nu$ is a weight of
$V_\lambda$ and $\nu\ne \lambda$.
It follows that $J_{\lambda,0}|_{\Sch_{\lambda}}
\cong \tilde Z_\lambda|_{\Sch_{\lambda}}$, and it suffices
to see that the sheaf in the right hand side of the latter isomorphism
has weight zero.
 This is clear from the construction of $\tilde Z_\lambda$
in \cite{zentr}. \epf

We now  define Weil perverse sheaves
 $\widetilde \nab^\lambda_{et} =\Phitil(\nab^{\lambda,\delta_\lambda})$,
 $\widetilde \del_\lambda^{et} =\Phitil(\del_{\lambda,\delta_\lambda})$,
 $\widetilde {IC}_\lambda=\Phit (E_{\lambda,\delta_\lambda})$ (the number $\delta_\la$
 was defined in \ref{nota}).
Lemma \ref{Phinab} and Corollary \ref{PhiIC}
show that these are Weil perverse sheaves
with underlying sheaves $\nab^\lambda$, $\del_\lambda$, $IC_\lambda$.

\begin{Prop}\label{Phidelchist}
The Weil sheaves  $\widetilde \nab^\lambda_{et}$,
 $\widetilde \del_\lambda^{et}$, $\ICtil_\lambda$
  are pure complexes of weight zero  on $\Fl^\lambda$
(the open stratum in their support).
\end{Prop}

{\em Proof.}\ It is enough to check the statement about
  $\widetilde \del_\lambda^{et}$, the rest of the claim follows
from this.
In view of \eqref{AvWakimot},
restriction of  $\widetilde \del_\lambda^{et}$ to the open stratum
$\Fl^\lambda$ in its support is isomorphic to the convolution
 $\widetilde\del_0^{et}* J_{\lambda, \delta_\lambda}$ restricted to
$\Fl^\lambda$.

 The $\LGO$ orbit $\Fl^\lambda$ fibers over $\LG/\LB$;
 Lemma  \ref{Waki_shriek} implies that the sheaf
$i_\lambda^*( J_{\lambda})$ is (up to shift and twist)
the pull-back of the standard sheaf attached to the $\LB$ orbit
in $\LG/\LB$ corresponding to
$w^\lambda\in W_f=\LB\backslash \LG/\LB$,
where $w^\lambda$ is the minimal length element
 such that $w^\lambda(\lambda)\in \Lambda^+$;
notice that $\ell(w^\lambda)=\delta_\lambda$.
In view of  Lemma \ref{Olachist} the claim follows now from the next Lemma.
 \epf

For $w\in W_f$ let $\i_w:(\LG/\LB)_w\imbed \LG/\LB$
 be the embedding of the corresponding $\LB$
orbit in $\LG/\LB$.

\begin{Lem}\label{weight_Whit}
For $w\in W_f$ we have an isomorphism of Weil sheaves
$$\widetilde{Wh}
*\i_{w!}(\cons[\ell(w)](\ell(w)/2))\cong \widetilde{Wh}(- \ell(w)/2).$$.
\end{Lem}

\proofpt We have $\i_{w_1!}(\cons[\ell(w_1)](\frac{\ell(w_1)}{2}))*
\i_{w_2!}(\cons[\ell(w_2)](\frac{\ell(w_2)}{2}))\cong
\i_{w_1w_2!}(\cons[\ell(w_1w_2)](\frac{\ell(w_1w_2)}{2}))$
 provided $\ell(w_1w_2)
=\ell(w_1)+\ell(w_2)$. Thus it is enough to prove the Lemma for $w=s$
 a simple reflection. In this case we have a short exact sequence of perverse
sheaves $$0\to L_s \to \i_{s!}(\cons[1](1/2))\to \delta_e(-1/2)\to 0$$
where $\delta_e$ is the skyscraper at the $\LB$-invariant point, and
$L_s$ is the constant sheaf (twisted and shifted) on the closure
of the orbit $\overline{(\LG/\LB)_s}\cong \Pone$.
It is easy to see that $Wh*L_s=0$, and $Wh*\delta_e\cong Wh$.
The Lemma follows. \epf

\medskip

We are now ready to express Ext's between coherent sheaves appearing
in the Positivity Lemma in terms of perverse sheaves.

Assume that $\iota_\la:x_\lambda\imbed \Fl^\la$ is defined over
$\Fq$. Then for a Weil perverse sheaf $\F$ on $\Fl$ the graded
vector spaces $\iota_\lambda^*(\F)$, $\iota_\lambda^!(\F)$ carry an
action of Frobenius. Let $\zeta_\lambda$ denote the constant by
which Frobenius acts on the one dimensional vector space
$\iota_\lambda^*(\widetilde \del^{et}
_{\lambda})[-d_\lambda](-\frac{d_\lambda}{2})$, where $d_\la=\dim
\Fl^\la$. Proposition \ref{Phidelchist} shows that $\zeta_\lambda$
is an algebraic number all of whose conjugates have absolute value
one (in fact, it is easy to see that $\zeta_\lambda$ is a root of
unity; we will not use this fact). Let
$\iota_\lambda^*(\F)_{[n]}^{(m)}$, $\iota_\lambda^!(\F)_{[n]}^{(m)}$
denote the $q^{m/2}\zeta_\lambda$
 eigenspace in the  $n$-th cohomology space  of the corresponding (co)stalk.

\begin{Cor}\label{graded_stalk}
We have canonical isomorphisms
$$\Hom^i_{D^{G\times \Gm}(\Nt)}(\Delta_{\lambda,n+\delta_\lambda},E_{\mu,\delta_\mu})
 = 
\iota_\lambda ^! (\widetilde {IC}_\mu) _{[i+d_\la]}^{(n+d_\la)};$$
$$\left[ \Hom^i_{D^{G\times \Gm}(\Nt)}(E_{\mu,\delta_\mu},\nab^{\lambda,\delta_\lambda+n})\right]^* =
 \iota_\lambda ^* (\widetilde {IC}_\mu)_{[-i-d_\la]}^{(n-d_\la)}
.$$

\end{Cor}

\proofpt The (co)stalk $\iota_\la^*(\F)$, $\iota_\la^!(\F)$ of an
object $\F\in D^b(\P_{asph})$ can be expressed in terms of Hom to
(from) a (co)standard sheaf, see \eqref{stalk_costalk}. Taking into
account the action of Frobenius, we see that if $\F\in
D^b(\P_{asph})$ is equipped with a Weil structure, then
$$ 
 \Hom^i\left(
\widetilde \del_{\lambda}^{et}
   ,  \F \left( \frac{n}{2} \right)\right) ^{Fr} \cong
   \iota_\lambda ^! (\F) _{[i+d_\la]}^{(n+d_\la)} ,$$
$$
 \left[ \Hom^{i} \left(\F,
 \widetilde
 \nab^{\lambda}_{et}(-\frac{n}{2})\right)^{Fr}\right]^*
\cong \iota^*_\la(\F)_{[-i-d_\la]}^{(n-d_\la)}
.$$

 We
 plug in $\F=\widetilde {IC}_\mu=\Phitil(E_{\mu,\delta_\mu})$;
 the claim follows then from Lemma \ref{gr_Fr}. \epf


\begin{Rem}
The content of the Corollary can be summarized as follows. The
graded components in the Ext spaces between irreducible exotic
sheaves and (co) standard sheaves can be expressed using weight
components of the (co)stalks of irreducible Iwahori-Whittaker
sheaves. Notice shifts by $\delta_\lambda$, $\delta_\mu$, which come
from the following three facts:

\begin{itemize}\item
a Wakimoto sheaf $J_\lambda$ is standard
(extension by zero from a cell) along the generic
fiber of the projection from its support to the affine Grassmannian
(Lemma \ref{Waki_shriek});
\item
  dimension of this generic fiber
equals $\delta_\lambda$ (see the proof
of Proposition \ref{Phidelchist});
\item
 the skyscraper sheaf at the zero-dimensional
Schubert cell in $\LG/\LB$ enters the Jordan-H\" older series of the
extension by zero of a weight zero sheaf on a  $d$-dimensional
Schubert cell with weight $-2d$ (Lemma \ref{weight_Whit}).
 \end{itemize}

This shift by $\delta_\lambda$ matches  the homological shift
by $\delta_\lambda$ appearing in section \ref{qua} (see Theorem \ref{ET} and
section \ref{33}); the latter is related to the homological
shift arising in the weak Borel-Weil
Theorem
(cf. the usual Borel-Weil Theorem, which implies that, over a field
of characteristic zero,
 $R\Gamma (G/B, \O(l\lambda))$ is concentrated in homological degree
$\delta_\lambda$ for  $l$ larger than the Coxeter number).

\end{Rem}

\subsection{Proof of the Positivity Lemma.} By Proposition
\ref{Phidelchist},
  $\widetilde {IC}_\mu|_{\Fl^\mu}$ is  pure  of
weight zero. Hence $\widetilde {IC}_\mu$ is pure of weight zero by
\cite{BBD},  Corollary 5.4.3. Thus
  the right hand side of
the first (respectively, second) isomorphism in Corollary
\ref{graded_stalk} vanishes for $i>n$ (respectively, $-i<n$)
 by
the definition of a pure complex (cf., e.g., \cite{BBD}, Corollary
5.1.9). \epf

\end{document}